\documentclass[11pt,draft]{article}
\usepackage{amsmath,amssymb,amsthm,enumerate,natbib,color}
\usepackage[latin1]{inputenc}
\usepackage{hyperref}
\def\Nset{{\mathbb{N}}}

\def\Rset{\mathbb R}

\newcommand{\ie}{{\em i.e.} }

\newcommand{\pscal}[2]{\langle #1, #2 \rangle}
\def\ind{{\bf 1}}
\newcommand{\Id}{\mathrm{Id}} % matrice identite

\newcommand{\gen}{\mathcal A}

\def\sfX{\mathsf{X}} % espace d'etat

\def\esp{\mathbb E}
\def\pr{\mathbb P}

\def\mcf{\mathcal{F}}

\newcounter{hypSM}

\newcounter{hyparf}

\newcounter{hypdr}

\newcounter{hypoconbis} \newcounter{saveconbis} \newcommand\debutA{
   \begin{list} {{\bf A\arabic{hypoconbis}}}
     {\usecounter{hypoconbis}}\setcounter{hypoconbis}{\value{saveconbis}}}
\newcommand\finA{
    \end{list}\setcounter{saveconbis}{\value{hypoconbis}}}

\newcounter{hypoconter} \newcounter{saveconter} \newcommand\debutB{
   \begin{list} {{\bf B\arabic{hypoconter}}}
     {\usecounter{hypoconter}}\setcounter{hypoconter}{\value{saveconter}}}
\newcommand\finB{
    \end{list}\setcounter{saveconter}{\value{hypoconter}}}

\newtheorem{theo}{Theorem}[section]
\newtheorem{lem}[theo]{Lemma}
\newtheorem{coro}[theo]{Corollary}
\newtheorem{prop}[theo]{Proposition}

 \theoremstyle{remark}

 \newcommand{\espCond}[2]{\esp \left[
    \left. #1\right| #2\right]}

\numberwithin{equation}{section}

\begin{document}
\bibliographystyle{plain} \title{Subgeometric rates of convergence of
  $f$-ergodic strong Markov processes} \author{Randal Douc \thanks{CMAP, Ecole
    Polytechnique, 91128 Palaiseau Cedex, France. douc@cmapx.polytechnique.fr}
  \and Gersende Fort \thanks{CNRS/LTCI, 46 rue Barrault, 75634 Paris Cedex 13,
    France. gfort@tsi.enst.fr} \and Arnaud Guillin \thanks{Ceremade, Université
    Paris Dauphine, Place Marechal de Lattre de Tassigny, 75775 Paris cedex 16,
    France. guillin@ceremade.dauphine.fr}} \date{} \maketitle

\begin{abstract}
 We provide a condition for $f$-ergodicity of strong
    Markov processes at a subgeometric rate. This condition is couched in terms
    of a supermartingale property for a functional of the Markov process.
    Equivalent formulations in terms of a drift inequality on the extended
    generator and on the resolvent kernel are given. Results related to
    $(f,r)$-regularity and to moderate deviation principle for integral
    (bounded) functional are also derived. Applications to specific processes
    are considered, including elliptic stochastic differential equation,
    Langevin diffusions, hypoelliptic stochastic damping Hamiltonian system and
    storage models.
\end{abstract}

\bigskip 

\noindent  {\bf Short title~:} Subgeometrically $f$-ergodic Markov processes.  

\noindent {\bf Keywords~:}
Subgeometric ergodicity, regularity, Foster's criterion, resolvent, moderate
deviations; Langevin diffusions, hypoelliptic diffusions, storage models.

\noindent {\bf MSC 2000 subject classifications~:} {\em Primary~:} 60J25, 37A25 {\em Secondary~:} 60F10,60J35, 60J60.

\clearpage
\newpage
 \section{Introduction}
 
 In the present paper, we study the recurrence of continuous-time Markov
 processes. More precisely, we provide a criterion that yields a precise
 control of a subgeometric moment of the return-time to a test-set.  The
 obtained result permits further quantitative analysis of characteristics such
 as the regularity of the process, the rate of convergence to the stationary
 state, and a moderate
 deviation principle. \\
 
 The stability and ergodic theory of continuous-time Markov processes has a
 large literature which is mainly devoted to the geometric case (also referred
 to as the exponential case). Meyn and Tweedie developed stability concepts for
 continuous-time Markov processes as well as simple criteria for
 non-explosivity, non-evanescence, Harris-recurrence, positive
 Harris-recurrence, ergodicity and geometric ergodicity
 ~\citep{Meyn:Tweedie:1993,Meyn:Tweedie:1993b,Meyn:Tweedie:1993c}.  Of
 particular importance in actually applying these concepts is the existence of
 verifiable conditions.  In the discrete-time context, development of
 Foster-Lyapunov type conditions on the transition kernel has provided such
 criteria (e.g.~\citep{MT93}).  In the continuous-time context, Foster-Lyapunov
 inequalities applied to the generator of the process play the same role.
 These criteria were successfully applied to the study of the solution to
 stochastic differential equations (see~\citep{Hasminskii:1980} and more
 recently,~\citep{Goldys:Maslowski:2006} and references therein).  Results
 relative to rates of convergence slower than geometric are not so well
 established.  In~\citep{Veretennikov:1997,Malyshkin:2001} (resp.
 \citep{Hou:Liu:Zhang:2005}), polynomial and sub-exponential ergodicity of
 stochastic differential equations (resp. sub-exponential ergodicity of queuing
 models) are addressed, but these results are quite model-specific. Fort and
 Roberts~\citep{Fort:Roberts:2005} are, to our best knowledge, the first to
 study the subgeometric ergodicity of general strong Markov processes.  Their
 conditions are in terms of subgeometric moment of the return-time to a
 test-set.  Fort and Roberts derive nested drift inequalities on the generator
 of the process that
 makes the result of practical interest in the polynomial case. \\
 
 One of the application of the condition we derive in the present paper makes
 the Fort-Roberts's theory applicable for more general subgeometric rates such
 as the logarithmic or the subexponential ones. It also provides criteria for
 the $(f,r)$-regularity of a process, a characteristic which is an extension of
 the regularity concept~\citep{Meyn:Tweedie:1993b}. We obtain theoretical
 results that are analogous to those in the discrete-time
 case~\citep{Tuominen:Tweedie:1994}. We then relate our condition to a
 criterion based on the generator of the process.  This criterion is the
 natural analogue of the Foster-Lyapunov condition for the geometric case; it
 also provides a single drift condition that generates the set of nested drift
 conditions by Fort-Roberts~\citep{Fort:Roberts:2005} for the polynomial case.
 Furthermore, it is analogous to the discrete-time version recently proposed by
 Douc-Fort-Moulines-Soulier \citep{Douc:Fort:Moulines:Soulier:2004}. \\
 In the literature, one approach for the theory of continuous-time Markov
 process is through the use of associated discrete-time chains~: the resolvent
 chains and/or a skeleton chain. We discuss how our condition is related to a
 subgeometric drift inequality for these discrete-time Markov chains. As a
 consequence, we state new limit theorems such as moderate deviations for
 integral of bounded functionals, thus weakening the conditions derived in
 Guillin-Wu~\citep{Guillin:2001,Wu:2001}.
 
 Our conditions are then successfully applied to various non trivial models:
 \textit{(a)} we first consider elliptic stochastic differential equations for
 which conditions on the drift function enable us to generalize results by
 Veretennikov \citep{Veretennikov:1997}, Ganidis-Roynette-Simonot
 \citep{Ganidis:Roynette:Simonot:1999} or Malyshkin \citep{Malyshkin:2001} (see
 also Pardoux-Veretennikov \citep{Pardoux:Veretennikov:2001} for a study of the
 regularity of the solution of the Poisson equation under this drift
 condition); \textit{(b)} we then study a "cold" Langevin tempered diffusion
 when the invariant target distribution is subexponential in the tails. This
 model is particularly useful in Markov Chain Monte Carlo method. Different
 regime of ergodicity (polynomial, subexponential or exponential) depending on
 the coldness of the diffusion term are exhibited, the different regimes are
 then characterized by the invariant target distribution. This study
 generalizes the Fort-Roberts' results, which consider the case when the
 target density is polynomial in the tails~\citep{Fort:Roberts:2005};
 \textit{(c)} we also give a toy hypoelliptic example, namely a stochastic
 damping Hamiltonian system, in the case when it cannot be geometrically
 ergodic. This model is shown to be polynomially ergodic (see
 Wu~\citep{Wu:2001} for the exponential case); \textit{(d)} we finally consider
 a simple compound Poisson-process driven Ornstein-Uhlenbeck process (relevant
 for recent studies in financial econometrics) with heavy tailed jump.  It is
 shown to be subgeometrically ergodic.

 Our approach may be considered as a probabilistic one. There are another ways
 to get subexponential rates of convergence (in total variation norm, in $L^2$
 or in entropy) such as those based on spectral techniques (as in
 \citep{Ganidis:Roynette:Simonot:1999}), or on functional inequalities (weak
 Poincar\'e inequalities~\citep{Rochner:Wang:2001} or weak logarithmic Sobolev
 inequalities~\citep{Cattiaux:Gentil:Guillin:2005}). These results are however
 not easy to compare to ours and we postpone a comparative utilization of these
 approaches to further research.

 Let us finally present the organization of the paper. Section 2 recalls basic
 definitions on Markov processes. The main results are given in
 Section~\ref{sec:MainResults}. All the proofs are postponed in appendix.
 Section~\ref{sec:Applications} is devoted to the examples and
 Section~\ref{discret} to a new moderate deviation principle.

\section{Definitions}
Let $(\Omega,\mcf,(\mcf_t)_{t\geq 0},(X_t)_{t\geq 0},(\pr_x)_{x \in \sfX})$ be  a Markov family on a locally compact and separable metric space $\sfX$ endowed with its Borel $\sigma$-field $\mathcal{B}(\sfX)$~: $(\Omega, \mcf)$ is  a measurable space, $(X_t)_{t \geq 0}$ is a  Markov process with respect to the filtration $(\mathcal{F}_t)_{t \geq 0}$ and $\pr_x$ (resp. $\esp_x$) denotes the canonical probability (resp. expectation) associated to the Markov process with initial distribution the point mass at $x$.  Throughout this paper, the process is assumed to be a  time-homogeneous   strong Markov process with cad-lag paths, and we denote by $(P_t)_{t \geq 0}$ the associated  transition function on $(\sfX, \mathcal{B}(\sfX))$. \\
Let $\Lambda_0$ denote the class of the measurable and nondecreasing functions
$r:[0, +\infty) \to [2,+\infty)$ such that $\log r(t) /t \downarrow 0$ as $t
\to +\infty$.  Let $\Lambda$ denote the class of positive measurable functions
$\bar r$, such that for some $r \in \Lambda_0$,
\[
0 < \liminf_t \frac{ \bar r(t)}{r(t)} \leq \limsup_t \frac{ \bar r(t)}{r(t)}  < \infty.
\]
$\Lambda$ is the class of the subgeometric rate functions and examples of
functions $\bar r \in \Lambda$ are
\[
\bar r(t) = t^\alpha \; (\log t)^\beta \; \exp(\gamma t^\delta)
\]
for $0< \delta <1$ and either $\gamma>0$, or $\gamma=0$ and $\alpha >0$, or
$\gamma=\alpha=0$ and $\beta \geq 0$. We are ultimately interested in
conditions implying that for all $x \in \sfX$
\begin{equation}
  \label{eq:DefiErg}
   \lim_{t \to +\infty} \ r(t) \ \ \| P^t(x,\cdot) - \pi(\cdot) \|_{f} = 0,
\end{equation}
where $r \in \Lambda$, $\pi$ is the (unique) invariant distribution of the process \ie\ $\pi P^t = \pi$ for all $t \geq 0$,  and for a signed measure $\mu$, $\| \mu \|_{f} = \sup_{|g| \leq f} |\mu(g)|$ where $f: \sfX \to [1, \infty)$ is a measurable function. When $f$ is the constant function $\ind$ ($\ind(t)=1$, $t\geq 0$), the $f$-norm is nothing more than the total variation norm. \\
To that goal, we will need different notions of regularity and stability of
continuous-time Markov processes and we briefly recall some basic definitions.  The
process is $\phi$-irreducible for some $\sigma$-finite measure $\phi$ on
$\mathcal{B}(\sfX)$ if $\phi(A)>0$ implies $\esp_x\left[\int_0^\infty
  \ind_A(X_s) \; ds \right]>0$ for all $x \in \sfX$. A $\phi$-irreducible
process possesses a maximal irreducibility measure $\psi$ such that $\phi$ is
absolutely continuous with respect to $ \psi$ for any other irreducibility
measure $\phi$~\citep{Nummelin:1984}. Maximal irreducibility measures are not
unique and are equivalent.  A set $A \in \mathcal{B}(\sfX)$ such that
$\psi(A)>0$ for some maximal irreducibility measure is said accessible; and
full if $\psi(A^c)=0$.  A measurable set $C$ is $\nu_a$-petite (or simply
petite) if there exist a probability measure $a$ on the Borel $\sigma$-field of
$[0, +\infty)$ and a non-trivial $\sigma$-finite measure $\nu_a$ on
$\mathcal{B}(\sfX)$ such that
  \[
  \forall x \in C, \qquad \qquad \int_0^{+\infty} P^t(x, \cdot) \; a(dt) \geq
  \nu_a(\cdot).
\]
For a $\psi$-irreducible process, an accessible closed petite set always
exists~\citep{Meyn:Tweedie:1993}.  A process is Harris-recurrent if, for some
$\sigma$-finite measure $\mu$, $\mu(A)>0$ implies that the event $\{
\int_0^\infty \ind_A(X_s) ds = \infty\}$ holds $\pr_x$-a.s. for all $x \in
\sfX$. Harris-recurrence trivially implies $\phi$-irreducibility.  A
Harris-recurrent right process possesses an invariant measure
$\pi$~\citep{Getoor:1980}; if $\pi$ is an invariant probability distribution, the
process is called positive Harris-recurrent.
A  $\phi$-irreducible  process is aperiodic if there exists an accessible $\nu_{\delta_m}$-petite set $C$ and $t_0$ such that for all $x \in C, t \geq t_0, P^t(x,C)>0$. A sufficient condition for a positive Harris-recurrent process to be  aperiodic is the existence of  some $\phi$-irreducible skeleton chain~\citep{Meyn:Tweedie:1993b}; recall that a skeleton  $P^m$ ($m>0$) is said $\phi$-irreducible if  there exists a $\sigma$-finite measure $\mu$ such that $\mu(A)>0$ implies $\forall x \in \sfX$, $\exists k \in \Nset$, $P^{km}(x, A) >0$~\citep{MT93}.  \\
A $\psi$-irreducible and aperiodic Markov process that verifies
(\ref{eq:DefiErg}) is said $f$-ergodic at a subgeometric rate (or simply
$f$-ergodic when $r=\ind$). When $r$ is of the form $r(t) = \kappa^t$ for some
$\kappa >1$, the process is said $f$-ergodic at a geometric rate.  In the
literature, criteria for the stability of Markov processes, when stability is
couched in terms of Harris-recurrence, positive Harris-recurrence,
$f$-ergodicity, with in this latter case, a mention of the rate of convergence,
are expressed in terms of hitting-times of some closed petite set. For any
$\delta>0$ and any closed set $C \in \mathcal{B}(\sfX)$, let
\[
\tau_C(\delta) = \inf \{t \geq \delta, X_t  \in C \},
\]
be the hitting-time on $C$ delayed by $\delta$ and define its $(f,r)$-modulated moment
\[
G_C(x,f,r; \delta) = \esp_x \left[ \int_0^{\tau_C(\delta)} \;  r(s) \; f(X_s) \; ds \right],
\]
where $f : \sfX \to [1, \infty)$ is a measurable function and $r: [0, +\infty)
\to (0, +\infty)$ is a rate function.  When $f = \ind$ (resp. $r= \ind$), this
moment is simply called the $r$-modulated (resp. $f$-modulated) moment.
Following discrete-time
usage~\citep{MT93,Tuominen:Tweedie:1994,Jarner:Roberts:2002}, we call a
measurable set $C$ $(f,r)$-regular if
\[
\sup_{x \in C} G_B(x,f,r; \delta) < \infty,
\]
for all $\delta >0$ and all accessible set $B$.    Criteria for Harris-recurrence and positive Harris-recurrence can be found in ~\citep[Theorems 1.1 and 1.2]{Meyn:Tweedie:1993}; ergodicity and  $f$-ergodicity are addressed in~\citep[Theorems 6.1 and 7.2]{Meyn:Tweedie:1993b}; criteria for geometric $f$-ergodicity at a geometric rate (resp. at a  subgeometric rate) are provided by~\citep[Theorem 7.4]{Down:Meyn:Tweedie:1995} (resp.~\citep[Theorem 1]{Fort:Roberts:2005}). A short review of these notions and results can be found in~\citep{Fort:Roberts:2005}.  \\
In many applications, these moments can not be explicitly calculated; a second
set of criteria based on the extended generator were thus derived for some of
the stability properties above. We postpone to Section~\ref{sec:Generator} a
review of the existing conditions.

%%%%%%%%%%%%%%%%%%%%%%%%%%%%%%%%%%%%%%%%%%%%%%%%%%%%%%%%%%%%%%
%%%%%%
%%%%%%   main results
%%%%%%
%%%%%%%%%%%%%%%%%%%%%%%%%%%%%%%%%%%%%%%%%%%%%%%%%%%%%%%%%%%%%%

\section{Main results}
\label{sec:MainResults}
Let us consider the following drift condition  towards  a closed petite set $C$.
 \begin{enumerate}
 \item []$\mathbf{D(C,V,\phi,b)}$: There exist a closed petite set $C$, a
   cad-lag function $V : \sfX \to [1,\infty)$, an increasing differentiable
   concave positive function $\phi: [1,\infty) \to (0,\infty)$ and a constant $b
   < \infty$ such that for any $s\geq 0$, $x \in \sfX$,
 \begin{equation}
   \label{eq:drift}
   \esp_x\left[V(X_s)\right]+\esp_x\left[ \int_0^s \phi \circ V(X_u)du\right]\leq V(x) +b \; \esp_x\left[ \int_0^s\ind_C(X_u)du\right]. 
 \end{equation}
 \end{enumerate}

Note that (\ref{eq:drift}) is equivalent to the condition that the functional
\[
s \mapsto V(X_s) -V(X_0) + \int_0^s \phi\circ V(X_u)du - b  \int_0^s\ind_C(X_u)du
\]
is, for all $x \in \sfX$, a $\pr_x$-supermartingale with respect to the
filtration $(\mathcal{F}_t)_{t \geq 0}$.

The main result of Section~\ref{sec:ModulatedMoments} is
Theorem~\ref{theo:ContMax} that states that this drift condition allows the
calculation of an upper bound for some $r$-modulated moment where $r \in
\Lambda$, and for some $f$-modulated moment, $f \geq 1$.  Using interpolating
inequalities, we obtain $(f,r)$-modulated moments for a wide family of pairs
$(f,r)$.  Section~\ref{sec:regularity} is devoted to $(f,r)$-regularity~: the
main result of this section is Proposition~\ref{prop:fr-regular} that
identifies $(f,r)$-regular sets from the condition $\mathbf{D(C,V,\phi,b)}$.
In Section~\ref{sec:ergodic}, we show that the drift condition
$\mathbf{D(C,V,\phi,b)}$ provides a simple sufficient condition for an
aperiodic strong Markov process to be $f$-ergodic at a subgeometric rate~: the
main result is Theorem~\ref{theo:FR_ergod} that builds on the work by Fort and
Roberts~\citep{Fort:Roberts:2005}. We provide in Section~\ref{sec:Generator} a
condition couched in terms of the extended generators that implies the drift
inequality $\mathbf{D(C,V,\phi,b)}$. This condition generalizes the condition
in \citep[Proposition 6]{Fort:Roberts:2005} that restricts to the polynomial
case, and reveals of great interest in many applications. We present in
Section~\ref{sec:equivResolv} the interplay between a drift condition on the
resolvent kernel and the drift condition $\mathbf{D(C,V,\phi,b)}$. \\
All the proofs are given in Appendix~\ref{app:Proofs}.

%%%%%%%%%%%%%%%%%%%%%%%%%%%%%
%%%modulated moments
%%%%%%%%%%%%%%%%%%%%%%%%%%%%%

\subsection{Modulated moments}
\label{sec:ModulatedMoments}
We show that $\mathbf{D(C,V,\phi,b)}$ is a simple condition that allows the
control of $f$-modulated moments and $r$-modulated moments, for a general rate
function $r \in \Lambda$, of the delayed hitting-time $\tau_C(\delta)$.  Let
\[
 H_\phi(u)=\int_1^u \frac{ds}{\phi(s)}, \qquad u \geq 1\ . 
\]
\begin{theo}
\label{theo:ContMax}
  Assume $\mathbf{D(C,V,\phi,b)}$.   
  \begin{enumerate}[i)]
\item \label{ContMax_Csq1} For all $x \in \sfX$ and  $\delta >0$,
$$
\esp_x \left[  \int_0 ^{\tau_C(\delta)} \phi \circ V(X_s) \; ds\right] \leq V(x)-1 + b \delta\ .
$$
  \item \label{ContMax_Csq2} For all $x \in \sfX$ and   $\delta >0$,
$$
\esp_x \left[ \int_0 ^{\tau_C(\delta)} \phi \circ H_\phi^{-1}(s)\;  ds \right] \leq V(x)-1+ \frac {b}{\phi(1)} \int_0 ^{\delta} \phi \circ H_\phi^{-1}(s)ds\ .
$$
\end{enumerate}
\end{theo}
%% The proof is given in Appendix~\ref{app:Proofs}.
The proof of Theorem~\ref{theo:ContMax} does not require $C$ to be
petite. Nevertheless, this petiteness property will be crucial in all the following results: we will see that this assumption  allows the extension of the above controls to those of modulated moments $\tau_B(\delta)$ for any accessible set $B$.  Theorem~\ref{theo:ContMax} gives the largest $f$-modulated and $r$-modulated
moments of $\tau_C(\delta)$ that can be deduced from $\mathbf{D(C,V,\phi,b)}$.
Interpolated $(f,r_f)$-modulated moments of $\tau_C(\delta)$ can easily be
obtained for a wide family of functions $1 \leq f \leq f_\ast$ (and,
equivalently, a wide family of rate functions $r(s) \leq r_\ast(s)$) where
\begin{equation}
  \label{eq:QuantMax}
  f_\ast = \phi\circ V, \qquad \qquad  r_\ast(s) = \phi \circ H_\phi^{-1}(s).
\end{equation}
 To that goal, we follow the same lines as in~\citep{Douc:Fort:Moulines:Soulier:2004} and \citep{Fort:Roberts:2005} and introduce  the pairs of Young's functions $(H_1, H_2)$  that, by definition,  satisfy the  property
\begin{equation}
  \label{eq:YoungIneq}
  x \; y \leq H_1(x) + H_2(y), \qquad \forall x,y \geq 0,
\end{equation}
and are invertible (see e.g \citep[Chapter 1]{Krasnoselskii:Rutickii:1961}). Let $\mathcal{I}$ be the
pairs of inverse Young's functions augmented with the pairs $(\Id, \ind)$ and
$(\ind, \Id)$. Examples of pairs $(H_1, H_2)$ are given in \citep{Douc:Fort:Moulines:Soulier:2004} and
\citep{Fort:Roberts:2005} while a general construction can be found in~\citep[Chapter
1]{Krasnoselskii:Rutickii:1961}. Corollary~\ref{coro:ContYoung} trivially results from
Theorem~\ref{theo:ContMax} and Eq. (\ref{eq:YoungIneq}).
\begin{coro}
  \label{coro:ContYoung}
  Assume $\mathbf{D(C,V,\phi,b)}$. For any pairs $(\Psi_1, \Psi_2) \in \mathcal{I}$ and all $\delta>0$,  
\[
\esp_x \left[  \int_0^{\tau_C(\delta)} \Psi_1\left(r_\ast(s)\right) \ \Psi_2\left(f_\ast(X_s) \right)  \; ds \right] \leq  2 (V(x)-1) +b \int_0^\delta \left( 1+ \frac{r_\ast(s)}{r_\ast(0)} \right) ds.
\]
\end{coro}
For two pairs $(\Psi_1, \Psi_2)$ and   $(\Psi_1', \Psi_2')$ in $\mathcal{I}$, if $\Psi_1(x) \leq \Psi_1'(x)$ for all large $x$, then $\Psi_2(y) \geq \Psi_2'(y)$ for all large $y$~\citep[Theorem 1.2.1]{Krasnoselskii:Rutickii:1961}. This shows that the rate  $\Psi_1\left(r_\ast(\cdot)\right)$ and the function $\Psi_2\left(f_\ast(\cdot)\right)$ have to be balanced~: the maximal rate function $r_\ast$ is associated to the  function $f$ with minimal growth in the range $1 \leq f \leq f_\ast$, that is with $f= \ind$; and the function with the largest rapidity of growth $f=f_\ast$ is associated to the minimal rate $r= \ind$.  \\
Theorem~\ref{theo:ContMax} and Corollary~\ref{coro:ContYoung} thus provides a
control of $(f,r)$-modulated moments; a simple condition for the rate $r$ to be
in the set $\Lambda$ of the subgeometric rate functions is recalled in the
following lemma~\citep[Lemmas 2.3 and 2.7]{Douc:Fort:Moulines:Soulier:2004}
\begin{lem}
 \label{lem:CS_rateLambda}
If $\lim_\infty \phi'=0$,  $r_\ast \in \Lambda$ and for all inverse Young function $\Psi_1$, $\Psi_1 \circ r_\ast \in \Lambda$.
\end{lem}

\begin{prop}
\label{prop:DandCaccessible}
Assume $\mathbf{D(C,V,\phi,b)}$. Then the process is $\psi$-irreducible.
If $\sup_C V< \infty$, 
\begin{enumerate}[(i)]
\item the level sets $\{ V \leq n\}$ are petite
and the union of these level sets is full.
\item there exists a closed accessible petite set $B$ such that
$\mathbf{D(B,V,\phi,b)}$ holds and $\sup_B V < \infty$. 
\end{enumerate}
\end{prop} 

As a consequence, when $\mathbf{D(C,V,\phi,b)}$ holds and $\sup_C V < \infty$,
we can assume without loss of generality that $C$ is accessible.

\subsection{$(f,r)$-regularity}
\label{sec:regularity}
Corollary~\ref{coro:ContYoung} shows that the drift condition
$\mathbf{D(C,V,\phi,b)}$ allows the control of modulated moments $G_C(x,f,r;
\delta)$, for all $\delta>0$ and a large family of pairs $(f,r)$.  Similar
modulated moments relative to any accessible set $B$ can be controlled provided
$\sup_{x \in C} G_C(x,f,r; \delta) < \infty$ for some $\delta>0$ (and thus any
$\delta >0$, as established in \citep[Lemma 20]{Fort:Roberts:2005}).
This naturally yields the notion of $(f,r)$-regular sets. The objective  of this section  is to identify regular sets from the drift condition $\mathbf{D(C,V,\phi,b)}$. \\

We start with a proposition that shows that the ``self-regularity'' of a closed
petite set $C$ actually implies $(f,r)$-regularity. This results extends
~\citep[Proposition 4.1]{Meyn:Tweedie:1993} (resp. \citep[Proposition 22]{Fort:Roberts:2005}) that
addresses the case $r = \ind$ (resp. $f= \ind$).  It also generalizes
~\citep[Proposition 23]{Fort:Roberts:2005} which concerns the case $r = \Psi_1 (r_\ast)$ and
$f = \Psi_2(f_\ast)$ for some pair $(\Psi_1, \Psi_2) \in \mathcal{I}$. This
proposition is the counterpart in the subexponential setting of the result by
Down-Meyn-Tweedie for the exponential case~\citep[Theorem 7.2]{Down:Meyn:Tweedie:1995}.

\begin{prop}
  \label{prop:FromCtoB}
  Let $f : \sfX \to [1, \infty)$ be a measurable function and $r \in \Lambda$
  be a subgeometric rate function.
  Assume that the process is $\psi$-irreducible and  $\sup_{x \in  C} G_C(x,f,r; \delta)<\infty$ for some (and thus any) $\delta >0$ and some closed petite set $C$.  \\
  For all accessible set $B \in \mathcal{B}(\sfX)$ and all $t\geq 0$, there
  exists a constant $c_{B,t}<\infty$ such that for all $x \in \sfX$,
 \begin{equation}
   \label{eq:fromCtoB}
   G_B(x,f,r; t)  \leq c_{B,t} \, G_C(x,f,r; \delta).
 \end{equation}
 Hence $C$ is $(f,r)$-regular.
\end{prop} 
\noindent 

\begin{prop}
  Assume that $\mathbf{D(C,V,\phi,b)}$ holds with $C,V, \phi$ such that
  $\sup_C V< \infty$ and $\lim_{+\infty} \phi'=0$. Then for any pair $(\Psi_1,
  \Psi_2) \in \mathcal{I}$, any accessible set $B$ and all $\delta >0$, there
  exists a finite constant $c$ such that
\[
\esp_x \left[ \int_0^{\tau_B(\delta)} \Psi_1\left(r_\ast(s)\right) \ 
  \Psi_2\left(f_\ast(X_s) \right) \; ds \right] \leq c \; V(x).
\]
Hence, any  $V$-level set $\{x \in \sfX, V(x) \leq v \}$ is $(f,r)$-regular for
all pairs $(f,r) = (\Psi_2 \circ f_\ast, \Psi_1 \circ r_\ast)$ with $(\Psi_1,
\Psi_2) \in \mathcal{I}$. 
\end{prop}
\begin{proof}
  By Corollary~\ref{coro:ContYoung}, $\sup_{x \in C} G_C(x,f,r; \delta)<\infty$
  for all $\delta>0$ provided the drift condition $\mathbf{D(C,V,\phi,b)}$
  holds and $\sup_C V< \infty$. Finally, $r = \Psi_1 \circ r_\ast$ for some
  inverse Young function $\Psi_1$ is a subgeometric rate if $\lim_{+\infty}
  \phi'=0$.  Proposition~\ref{prop:FromCtoB} thus implies that the level sets
  of $V$ are $(f,r)$-regular sets.
\end{proof}
We now establish a general result that extends to continuous-time Markov
processes, part of~\citep[Theorem 2.1]{Tuominen:Tweedie:1994} relative to discrete-time Markov
chain.  In the case $r = \ind$, some of these equivalences are proved
in~\citep{Meyn:Tweedie:1993} for continuous-time strong Markov processes.

\begin{prop}
  \label{prop:fr-regular}
  Let $f : \sfX \to [1, \infty)$ be a measurable function and $r \in \Lambda$
  be a subgeometric rate function.  Assume that the process is
  $\psi$-irreducible. The following conditions are equivalent
\begin{enumerate}[i)]
\item There exist a closed petite set $C$  and $\delta >0$ such that $\sup_C G_C(x,f,r; \delta)< \infty$.
\item There exists a $(f,r)$-regular closed set which is accessible.
\item There exists a full set $\mathcal{S}_\Psi$ which is the union of a
  countable number of $(f,r)$-regular sets.
\end{enumerate}
\end{prop}
Theorem~\ref{theo:ContMax} proves that these equivalent conditions are verified
provided  $\mathbf{D(C,V,\phi,b)}$ holds,  $\sup_C V< \infty$ and $\lim_\infty \phi'=0$.\\
We conclude this section by establishing that under mild additional conditions,
the drift condition $\mathbf{D}$ also yields controls of modulated moments for
the skeleton chains. For all $m>0$, let $T_{m,C}$ be the return-time to $C$ of
the skeleton chain $P^m$,
$$
T_{m,C}= \inf\{k \geq 1,  X_{mk} \in C\}.
$$

\begin{prop}
\label{prop:lienContDiscret} 
Assume that $\mathbf{D(C,V,\phi,b)}$ holds with $\sup_C V< \infty$, and some
skeleton chain is irreducible.  For all $m >0$ and any accessible set $B$,
there exist constants $c_i, 1 \leq i \leq 4$, such that for all $x \in \sfX$,
\[
\esp_x\left[ \sum_{k=0}^{T_{m,B}-1} \phi \circ V(X_{mk})\right] \leq c_1 \;
\esp_x \left[ \int_0^{T_{m,B}} \phi \circ V(X_{sm}) \; ds \right] \leq c_2 \;
V(x),
\]
and
$$
\forall x \in \sfX,\quad  \ \esp_x\left[ \sum_{k=0}^{T_{m,B}-1} r_\ast(km) \right] \leq   c_3 \ \esp_x \left[ \int_0^{\tau_B(\delta)}  \; r_\ast(s) \; ds \right] \leq  c_4 \; V(x).
$$
\end{prop}
We will see in the last section that this proposition which clearly links the
behavior of the skeleton chain to that of the initial process leads to new
limit theorems such as moderate deviations.  It will also imply interesting
applications to averaging principle.

\subsection{ $f$-ergodicity at a subgeometric rate} 
\label{sec:ergodic}
From the control of $x \mapsto G_C(x,f,r; \delta)$ where $C$ is a  closed petite set, we are able to deduce results on the ergodic behavior of the strong Markov process. \\
The first result concerns the existence of an invariant probability
distribution $\pi$ and shows that the drift condition $\mathbf{D(C,V,\phi,b)}$
provides a simple tool when identifying the set of the $\pi$-integrable
functions. The second one states that the Markov process converges in $f$-norm
to the invariant probability measure $\pi$, for a wide family of functions $1
\leq f \leq f_\ast$ and a wide family of rate functions $ r_f \leq
r_\ast$.

\begin{prop}
  \label{prop:CS_PiInteg}
  Assume $\mathbf{D(C,V,\phi,b)}$ and $\sup_C V< \infty$. Then the process is
  positive Harris-recurrent with an invariant probability measure $\pi$ such
  that $\pi(\phi \circ V) < \infty$.
\end{prop}
Proposition~\ref{prop:CS_PiInteg} results from \citep[Theorems 1.1 and
1.2]{Meyn:Tweedie:1993} and Theorem~\ref{theo:ContMax}(\ref{ContMax_Csq1}).
It is known that positive Harris-recurrence does not necessarily imply ergodicity and  aperiodicity is required~\citep[Proposition 6.1]{Meyn:Tweedie:1993b}; similar conditions are required in the discrete-time case~\citep{MT93}.  In the present case, we have more information than  positive Harris-recurrence and thus, we are able to establish  $f$-ergodicity at a subgeometric rate.    \\
For a sequence $r \in \Lambda$, define $r^0(t) = \int_0^t r(s) \ ds$, and, for
a differentiable rate function $r$, set $\partial r(t) = \frac{dr(t)}{dt}$.

\begin{theo}
\label{theo:FR_ergod}
  Assume that
  \begin{enumerate}[(i)]
  \item some skeleton chain  is irreducible.
  \item the condition $\mathbf{D(C,V,\phi,b)}$ holds with $C, V, \phi$ such that $\sup_C V<\infty$  and $\lim_{+\infty} \phi'=0$.
  \end{enumerate}
  For any pair $\Psi=(\Psi_1, \Psi_2) \in \mathcal{I}$ and any probability
  distribution $\lambda$ satisfying $\lambda(V) <\infty$,
\begin{equation}
  \label{eq:Ergodicite}
  \lim_{t \to +\infty} \ \ \left\{ \Psi_1(r_\ast(t)) \vee 1 \right\} \ \  \int_\sfX \; \lambda(dx) \; \|
P^t(x, \cdot) - \pi(\cdot) \|_{\Psi_2(f_\ast) \vee 1 } = 0,
\end{equation}
where $r_\ast$ and $f_\ast$ are given by (\ref{eq:QuantMax}) and $\mathcal{I}$
is defined in Section~\ref{sec:ModulatedMoments}. Furthermore, there exist
finite constants $C_{\Psi,i}$ such that for all $t \geq 0$ and all $x \in
\sfX$,
\begin{eqnarray}
  \label{eq:ControleExplicite1}
 && \left\{ \Psi_1(r_\ast(t)) \vee 1 \right\} \ \ \|
P^t(x, \cdot) - \pi(\cdot) \|_{\Psi_2(f_\ast) \vee 1 } \leq C_{\Psi,1} \; V(x), \\   
\label{eq:ControleExplicite3}
&&\int_0^\infty  \left\{  \Psi_1(r_\ast(t)) \vee 1  \right\} \ \ \ \|
P^t(x, \cdot) - P^t(y, \cdot)  \|_{\Psi_2(f_\ast) \vee 1 } \; dt  \leq
C_{\Psi,2} \; \{ V(x) + V(y) \};
\end{eqnarray}
and if $\partial[\Psi_1(r_\ast)] \in \Lambda$, there exists a finite constant $C_{\Psi,3}$ such that
for all $t \geq 0$,
\begin{equation}
  \label{eq:ControleExplicite2}
\int_0^\infty \left\{ \partial[\Psi_1(r_\ast)](t) \vee 1  \right\} \ \ \ \| 
P^t(x, \cdot) - \pi(\cdot) \|_{\Psi_2(f_\ast) \vee 1 } \; dt \leq C_{\Psi,3} \; V(x).
\end{equation}
\end{theo}
\noindent The limit~(\ref{eq:Ergodicite}) is a direct application of~\citep[Theorem 1]{Fort:Roberts:2005} while (\ref{eq:ControleExplicite1}) to (\ref{eq:ControleExplicite2}) are, to our best knowledge, new results. The proof of this theorem is detailed in Appendix~\ref{app:Proofs}. \\
As already commented in \citep{Fort:Roberts:2005}, Eq.~(\ref{eq:Ergodicite}) shows that the rate of convergence and the norm in which convergence occurs have to be balanced~: if $\Psi_1$ strongly increases at infinity  then $\Psi_2$ slowly  increases (see~\citep{Krasnoselskii:Rutickii:1961} and the comments in Section~\ref{sec:ModulatedMoments}). Hence, the stronger the norm, the weaker the rate and conversely. The maximal rate of convergence is achieved with the total variation norm ($\Psi_2 \circ f_\ast = \ind$) and the minimal one ($\Psi_2 \circ r_\ast = \ind$) is achieved  with the  $f_\ast$-norm. Hence, the drift condition $\mathbf{D(C,V,\phi,b)}$ directly provides two major informations: the largest rate of convergence $r_\ast = \phi \circ H_\phi^{-1}$ is given by the concave function $\phi$ and the largest norm of convergence $\|\cdot \|_{f_\ast}$  is given by the pair $(\phi,V)$. \\
Eqs. (\ref{eq:ControleExplicite1}) to (\ref{eq:ControleExplicite2}) are, to our best knowledge, the first results that address the dependence upon the initial point in the ergodic behavior. When applied to discrete-time Markov chains, (\ref{eq:ControleExplicite1}) to (\ref{eq:ControleExplicite2}) coincide with resp.~\citep[Theorems 2.1, 4.1, 4.2]{Tuominen:Tweedie:1994} (the dependence upon $x$ can be read from the proof of these theorems; the details are also provided in~\citep[Chapter 3]{Fort:2001}). These results for the discrete-time case and the definition of the set $\mathcal{S}_\psi$ in~\citep[Theorem 1]{Fort:Roberts:2005}  suggest that in (\ref{eq:ControleExplicite1}), the minimal dependence in the starting value $x$ is of the form $G_C(x,\Psi_2(f_\ast), \Psi_1(r_\ast); \delta)$. Similar expressions can be predicted  for (\ref{eq:ControleExplicite3}) and (\ref{eq:ControleExplicite2}). The proof of this assertion and the explicit construction of the constants $C_{\Psi,i}$  in terms of the quantities appearing in the assumptions are beyond the scope of this paper. Currently in progress is work on explicit control of subgeometric ergodicity for strong Markov processes. \\
In the examples given in Section~\ref{sec:Applications}, we will see that the
pair $(\phi,V)$ that solves $\mathbf{D(C,V,\phi,b)}$ is not unique. Roughly
speaking, we read from Theorem~\ref{theo:FR_ergod} that $\phi$ is related to
the rate of convergence in total variation norm, while $V$ is the dependence
upon the initial point in the control of convergence. As a consequence, the
rate of convergence $r_\phi(t)$ and the dependence $V(x)$ can be balanced to
make the bounds (\ref{eq:ControleExplicite1}) to (\ref{eq:ControleExplicite2})
minimal.  In Section~\ref{sec:Applications}, we will give some examples (on
$\sfX = \Rset^n$), where both a pair of polynomially increasing functions and a
pair of subgeometrically increasing functions can be found.  One then
immediately remarks that the stronger the control in the initial point is, the
stronger the decay in time is for a given norm. It stresses once again the
interest for exact constant in our controls to decide which "ergodicity" to use
to reach a certain level. The fact that the pair ($\phi, V$) is not unique
shows that the drift condition only provides an upper bound of the true rate of
convergence.  Nevertheless, in many applications, we are able to prove that the
true rate belongs to the exhibited class of rate functions (see for example,
section~\ref{sec:Ex2}).

%%%%%%%%%%%%%%%%%%%%
%%generator
%%%%%%%%%%%%%%%%%%%%

\subsection{Generator and drift inequality (\ref{eq:drift})}
\label{sec:Generator}
The drift condition $\mathbf{D(C,V,\phi,b)}$ may not be easy to derive since it
is couched in terms of the process itself. The main goal of this section is to
provide an easier path to ensure subgeometric ergodicity, which is moreover the
usual form of conditions adopted on earlier paper to address different classes
of stability.  Namely we will use the formalism of the extended
generator~\citep[Def. 1.15.15]{Davis:1993}.

Let $\mathcal{D}(\gen)$ denote the set of measurable functions $f: \sfX \to
\Rset$ with the following property: there exists a measurable function $h :
\sfX \to \Rset$ such that the function $t \mapsto h(X_t)$ is integrable
$\pr_x$-a.s. for each $x \in \sfX$ and the process
\begin{equation}
  \label{eq:DefiCtF}
  t \mapsto f(X_t) -f(X_0) - \int_0^t  h(X_s) ds
\end{equation}
is a $\pr_x$-local martingale for all $x$. Then we write $h = \gen f$, and $f$ is said in the domain of the extended generator $(\gen,
\mathcal{D}(\gen))$ of the process $X$. The condition
(\ref{eq:drift}) looks like a Dynkin formula. This is the reason why we want it
to hold as widely as possible, thus justifying the interest in the extended
generator concept.
\begin{theo}
\label{theo:generator2df}
Assume that there exist a closed petite set $C$, a cad-lag function $V:
\sfX\to[1,\infty)$ with $V \in \mathcal{D}(\gen)$, an increasing differentiable
concave positive function $\phi:[1,\infty)\to(0,\infty)$ and a constant $b <
\infty$ such that for all $x \in \sfX$,
\begin{equation}\label{eq:gendrift}
\gen V(x)\le -\phi\circ V(x)+b\ind_C(x). 
\end{equation}
Then $\mathbf{D(C,V,\phi,b)}$ holds.
\end{theo}
The proof is in Section \ref{proof:generator}. The extended generator is less
restrictive than the infinitesimal generator $\tilde \gen$~: if $f$ is in the
domain of $\tilde \gen$, then the process (\ref{eq:DefiCtF}) is a martingale
and $f$ is in the domain of $\gen$ (see e.g.  \citep[Proposition
1.14.13]{Davis:1993}). In particular, it is often quite difficult to characterize
the domain of $\tilde \gen$ but there may be (and are, in the applications of
Section~\ref{sec:Applications}) easily checked
sufficient conditions for membership of $\mathcal{D}(\gen)$. \\
This drift condition naturally inserts in the existing literature, that
addresses criteria for non-explosivity, recurrence, polynomial ergodicity,
geometric and uniform ergodicity.  More precisely, Meyn and Tweedie provide
conditions for non-explosion, recurrence, positive-Harris recurrence and
$V$-ergodicity at a subgeometric rate, respectively of the form
\begin{eqnarray}
  \gen V(x) &\leq& c V(x), \label{eq:NonExplosion} \\
\gen V(x)  &\leq & c \ind_C(x),  \label{eq:Recurrence} \\
\gen V(x)  &\leq & -c f(x) + b \ind_C(x), \label{eq:PositiveHR} \\
\gen V(x)  &\leq & -c V(x) + b \ind_C(x) \label{eq:GeometricErgo}
\end{eqnarray}
for some positive constants $b,c<\infty$ and a measurable function $f\geq 1$
(see~\citep[Conditions (CD0) to (CD3)]{Meyn:Tweedie:1993c}; see also \citep{Down:Meyn:Tweedie:1995} for the
condition~(\ref{eq:GeometricErgo})). These criteria are similar to some
conditions provided by~\citep{Hasminskii:1980} for the stability of stochastic
differential equations.  The drift inequality (\ref{eq:GeometricErgo}) is the
limit of our approach, since it corresponds to
(\ref{eq:gendrift}) with $\phi(v)  \propto v$. \\
In a recent work, Fort and Roberts \citep{Fort:Roberts:2005} considered a family of drift
condition that implies $f$-ergodicity at a polynomial rate~: namely, there
exist $0 < \alpha<1$, $b>0$ such that for all $\alpha\le\eta\le1$, there exists
$c_\eta>0$ such that
\begin{equation}
  \label{eq:GenDriftFR}
  \gen V^\eta(x) \le -c_\eta V^{\eta-\alpha}(x)+b 1_C(x).
\end{equation}
Our drift condition (\ref{eq:gendrift}) with $\phi (v) \propto v^{1-\alpha}$
yields the same results as those provided in \citep[Theorem 1, Lemma 25,
Proposition 26]{Fort:Roberts:2005} (see Theorem~\ref{theo:FR_ergod} and
Proposition~\ref{theo:ContMax}). Hence, the drift inequality
(\ref{eq:gendrift}) that addresses subgeometric ergodicity generalizes the
criterion for polynomial ergodicity proposed by Fort-Roberts.  The comparison
of the Fort-Roberts nested drift conditions (\ref{eq:GenDriftFR}) and our
single drift condition can be more explicit when $V \in \mathcal{D}(\gen)$ and
the process (\ref{eq:DefiCtF}) is a $\pr_x$-martingale for all $x$. In that
case, it is easily seen that the single drift condition implies the nested
drift conditions. The martingale property is equivalent to
\[
t \mapsto \exp \left( \ln V(X_t) - \ln V(X_0) -\int_0^t H(\ln V)(X_s) ds\right)
\]
is a $\pr_x$-martingale for all $x$, where $H(\ln V) = V^{-1} \gen
V$~\citep{Feng:Kurtz:2006}.  Furthermore, $H(\ln V) \leq -g +s$ if and only if
\[
t \mapsto \exp \left( \ln V(X_t) - \ln V(X_0) - \int_0^t \{-g(X_u) +s(X_u)\}
  du\right)
\]
is a $\pr_x$-supermartingale for all $x$~\citep{Feng:Kurtz:2006}. As a
consequence, if $V^\eta$ is in the domain of $\gen$ for all $0 \leq \eta \leq
1$ then the Jensen's inequality yields $H(\eta \ln V) \leq \eta \exp(-\alpha
\ln V) + b \eta \exp(- \ln V) 1_C$ which in turn implies (\ref{eq:GenDriftFR}).

\subsection{Resolvent and drift inequality (\ref{eq:drift})}
\label{sec:equivResolv}
One of the approaches for studying the stability and ergodic theory of
continuous time Markov processes consists in making use of the associated
discrete time resolvent chains. This allows to take profit of the analysis of
discrete time Markov chains which is quite well understood
(\citep{Nummelin:1984,MT93}) and then to transfer properties established in terms of the
resolvent or ``generalised resolvent'' kernel (see for e.g.
\citep{Meyn:Tweedie:1993}) to the Markov process itself. Following the
discussion (done for exponentially ergodic Markov process) by
Down-Meyn-Tweedie~\citep[Th.5.1]{Down:Meyn:Tweedie:1995} and extending it to the
subgeometric case, we will now link the drift condition
$\mathbf{D(C,V,\phi,b,\beta)}$ associated to the Markov process to a drift
condition associated to the discrete time resolvent chain.

More precisely, define, for $\beta >0$, the resolvent kernel $R_\beta$ by $
R_\beta(x,A)=\int_0^{\infty} \beta e^{-\beta t}P^t(x,A)dt $ and consider the
following drift condition associated to the resolvent kernel.
\begin{enumerate}
\item []$\mathbf{\check D(C,V,\phi,b,\beta)}$: There exist a petite set
  $C$, a function $V : \sfX \to [1,\infty)$, an increasing
  differentiable concave positive function $\phi: [1,\infty) \to (0,\infty)$
  and a constant $b < \infty$ such that for any $x \in \sfX$,
\begin{equation}
  R_\beta V(x) \leq V(x)- \phi \circ V(x)+b \ind_C(x). 
\end{equation}
\end{enumerate}

The following result ensures that drift conditions expressed in terms of the
resolvent kernel or of the Markov process are essentially equivalent. This
theorem parallels Theorem 5.1. by Down-Meyn-Tweedie \citep{Down:Meyn:Tweedie:1995} for
exponentially ergodic Markov processes.
\begin{theo}
\label{theo:equivResolv}
  \begin{itemize}
  \item [(i)] Assume $\mathbf{\check D(C,V, \phi,b,\beta)}$ where $C$ is a
    closed set and $R_\beta V$ is a cad-lag function. Then $\mathbf{D(C,R_\beta
      V, \beta \phi,\beta b)}$ holds.
  \item [(ii)] Assume $\mathbf{D(C,V, \phi,b)}$ with $\sup_C V <\infty$.  Then,
    for all $\epsilon>0$, there exists a constant $c$ such that for all $x \in
    \sfX$,
$$
W(x) \leq (1+\epsilon)V(x)+ c \quad \mbox{and} \quad \lim_{t \to \infty} \frac{r_{\check \phi}(t)}{ r_{\phi}((1+\epsilon) t)}=1+\epsilon
$$
such that $\mathbf{\check D(\check C,W, \check \phi,\check b,\beta)}$ holds.
  \end{itemize}
 \end{theo}

The proof is given in Section \ref{proof:resolv}.

%%%%%%%%%%%%%%%%%%%%%%%%%%%%%%%%%%%%%%
%%%
%%% exemples
%%%
%%%%%%%%%%%%%%%%%%%%%%%%%%%%%%%%%%%%%%

\section{Examples}
\label{sec:Applications}
In this section, $\sfX = \Rset^n$. Vectors are intended as column vectors,
$|x|$ and $\pscal{\cdot}{\cdot}$ denote respectively the Euclidean norm and the
scalar product.  For a matrix $a$, $|a| = \left(\sum_{i,j} a_{i,j}^2 \right)^{1/2}$, $\mathrm{Tr}(a)$ stands for the
trace of the matrix and $a'$ the matrix transpose.  $\Id_n$ is the $n \times n$ identity matrix. If $V$ is a twice
continuously differentiable function with respect to $x \in \Rset^n$, $\partial
V$ (or $\partial_x V$ when confusion is possible) denotes its gradient, and
$\partial^2 V$ its Hessian. \\ For a set $A$, $A^c$ is its complement in
$\Rset^n$.

Four applications are considered: we first analyze general elliptic diffusions on
$\Rset^n$ such that the drift coefficient verifies a contraction condition of
the form $\pscal{b(x)}{x} \leq -r |x|^{1-p}$ for all large $x$, where $0<p<1$.
We then consider a Langevin diffusion on $\Rset^n$ having an  invariant
distribution which is super-exponential in the tails, and show that the rate of
convergence can be modified by ``heating'' the diffusion. The method is however not limited to elliptic diffusions but can also be of use in the hypoelliptic case: we consider as an illustration a simple stochastic damping Hamiltonian system which cannot be exponentially ergodic but is shown to be subexponentially ergodic. We finally study a
compound Poisson-process driven Ornstein-Uhlenbeck process when the
distribution of the jump is heavy tailed.

Queuing theory is another important field of application for our theory.
We do not discuss here this field of applications. %  for our theory which is queuing theory that will be, to save space and to the very nature of the verification of our hypothesis in this case, treated in another paper.
This will be done in a forthcoming paper, which will also include a comparison
of our results to those by~\citep{Dai:Meyn:1995,Hou:Liu:Zhang:2005}. Techniques in Dai-Meyn
\citep{Dai:Meyn:1995} differ from ours since they are based on fluid limits.
Concerning~\citep{Hou:Liu:Zhang:2005}, our conditions are more general; indeed the authors
assume that there exists a state $x_0$ such that whenever the Markov process
hits $x_0$, it will sojourn there for a random time that is positive with
probability $1$, \citep[Assumption 1.1]{Hou:Liu:Zhang:2005}. This assumption makes their
results unavailable for the applications we now consider.

%%%%%%%%%%%%%%%%
%% diffusions
%%%%%%%%%%%%%%%

\subsection{Elliptic diffusions on $\Rset^n$}
\label{sec:Ex1}
Consider the stochastic integral equation of the form
\begin{equation}
  \label{eq:StochEquDiff}
  X_t = X_0 + \int_0^t b(X_s) ds + \int_0^t \sigma(X_s) dB_s,
\end{equation}
where $X_t \in \Rset^n$, $b : \Rset^n \to \Rset^n$ and $\sigma : \Rset^n \to \Rset^{n \times n}$  are measurable functions, and $\{B_t\}_t$ is a $n$-dimensional Brownian motion. Assume that $b: \Rset^n \to \Rset^n$ and $\sigma:  \Rset^n \to \Rset^{n \times n}$ are functions satisfying
\debutA
\item \label{Ex1:A1}
$\sigma$ is bounded and  $b$ and $\sigma$ are  locally Lipschitz~: for any $l> 0$, there exists a finite constant $c_l$ such that for all $|x| \leq l, |y| \leq l$,
    \begin{equation}
      \label{eq:Ex1_Lips}
      | b(x) - b(y) | + |\sigma(x) - \sigma(y)| \leq c_l |x-y|.
    \end{equation}
\finA
Let $a(x ) =\sigma(x) \sigma(x)'$ be the diffusion matrix. We assume that
\debutA
\item \label{Ex1:A2}  
  \begin{enumerate}[(i)]
  \item \label{Ex1:A2i} $a(x)$ is non-singular~: the smallest eigenvalue of the
    diffusion matrix $a(x)$ is bounded away from zero in every bounded domain.
\item \label{Ex1:A2iii} there exist $0<p<1$, $r>0$ and $M$ such that for all $|x| \geq M$, $\pscal{b(x)}{x} \leq -r |x|^{1-p}$.
  \end{enumerate}
  \finA Note that under A\ref{Ex1:A1}, $\Lambda = n^{-1} \sup_{x \in \Rset^n}
  \mathrm{Tr}(a(x))$ and $\lambda_+ = \sup_{x \neq 0} \pscal{a(x)
    \frac{x}{|x|}}{\frac{x}{|x|}}$ are finite. Moreover, since under
  A\ref{Ex1:A1} $\sigma$ is continuous, the assumption
  A\ref{Ex1:A2}(\ref{Ex1:A2i}) is equivalent to the condition
  $\mathrm{det}(\sigma(x)) \neq 0$ for all $x$.

Under A\ref{Ex1:A1}, it is possible to define continuous functions $b_l$ and $\sigma_l$ that satisfy the at most linear increasing
\[
| b_l(x)| + |\sigma_l(x)| \leq c_l (1+|x|), \qquad \forall x \in \Rset^n,
\]
the Lipschitz condition (\ref{eq:Ex1_Lips}) on the whole state space, and are
such that $b_l = b$ and $\sigma_l = \sigma$ on the cylinder $\{x \in \Rset^n,
|x| < l\}$. The stochastic equation (\ref{eq:StochEquDiff}) has a unique
$t$-continuous solution $\{X_t^{(l)} \}_t$, when $b$ and $\sigma$ are replaced
by $b_l$ and $\sigma_l$~\citep[Theorem 3.3.2]{Hasminskii:1980}. The first exit times of
$\{X_t^{(m)} \}_t$ from $\{x \in \Rset^n, |x| < l\}$ are identical for all $m
\geq l$ (and is thus denoted $\zeta_l$). This allows the construction of a
process $\{X_t\}_t$ that satisfies (\ref{eq:StochEquDiff}) up to the explosion
time $\zeta= \lim_l \zeta_l$. If $\zeta=+\infty$ a.s., $\{X_t\}_t$ is a.s.
defined for all $t\geq 0$ and the process is said regular. Under the stated
assumptions, an easy to check sufficient condition for regularity relies on the
operator $L$ that acts on function $V : \Rset^n$, $x
\mapsto V(x)$ that are 
twice continuously differentiable with respect to $x$:
\begin{equation}
  \label{eq:Ex1:Loperator}
  LV(x) =  \pscal{b(x)}{\partial V(x)} + \frac{1}{2} \mathrm{Tr}\left(a(x) \ \  \partial^2 V(x) \right).
\end{equation}
The process is regular if there exists a non-negative twice-continuously
differentiable function $V$ on $\Rset^n$ such that for some finite $c$, $LV
\leq c V$ on $\Rset^n$ and $\inf_{|x| >R} V(x) \to \infty$ as $R \to
\infty$~\citep[Theorem 3.4.1.]{Hasminskii:1980}.  Under A\ref{Ex1:A2}(\ref{Ex1:A2iii}), it
is trivial to verify that by setting $V(x) = 1+ |x|^2$,
\begin{equation}
  \label{eq:Ex1_LVpower2}
  LV(x) \leq \left\{
  \begin{array}{ll}
-2 r |x|^{1-p} + n \Lambda, \quad & \text{if $|x| \geq M$}, \\
2 M \ \ \sup_{|x| \leq M} |b(x)| + n \Lambda \quad &\text{otherwise}.
  \end{array} \right.
\end{equation}
This shows that the process is regular. Consequently, there exists a solution
to (\ref{eq:StochEquDiff}), which is an almost surely continuous stochastic
process and is unique up to equivalence. This solution is an homogeneous Markov
process whose transition functions are Feller functions~\citep[Theorem
3.4.1]{Hasminskii:1980}. Hence, it is strongly Markovian, as a right-continuous Markov
process with Feller transition functions.  We now discuss the existence of an
irreducible skeleton $P^m$ and the petiteness property of the compact sets.
All of these properties deduce from the existence of an unique invariant
probability distribution $\pi$.

\begin{prop}
  \label{prop:Ex1_ResPrelim}
 Under A\ref{Ex1:A1}-A\ref{Ex1:A2}, $X$ possesses an unique invariant
  probability measure $\pi$.  $\pi$ is a maximal irreducibility measure and any
  skeleton $P^m$ is irreducible. Furthermore, the compact sets are closed
  petite sets.
\end{prop}
\begin{proof}
  By (\ref{eq:Ex1_LVpower2}), \citep[Theorem 3.7.1]{Hasminskii:1980} and its corollary
  2~\citep[p. 99]{Hasminskii:1980}, there exists a bounded domain $U$ with regular
  boundary and a finite constant $c$ such that for all $x \in U^c$,
  $\esp_x[T_{U^c}] < \infty$ and for any compact $\mathcal{K} \subset \Rset^n$,
  $\sup_{x \in \mathcal{K}} \esp_x \left[T_{U^c} \right] < \infty$, where
\[
T_U = \inf \{ t \geq 0, X_t \notin U \}.
\]
Since the diffusion matrix $a(x)$ is non-singular, we deduce from~\citep[Theorem
4.4.1 and Corollary 2 p.123]{Hasminskii:1980} that the process possesses an unique invariant
probability distribution $\pi$.  \citep[Lemma 4.6.5]{Hasminskii:1980} implies that any
skeleton is $\phi$-irreducible, with an irreducibility measure absolutely
continuous with respect to the Lebesgue measure. By \citep[Lemma 4.6.1]{Hasminskii:1980}, the support of $\pi$ has non-empty interior;
since the process is $\psi$-irreducible and has the Feller property, all
compact subsets of $\Rset^n$ are petite (this assertion can be proved in
exactly the same way as in the discrete-parameter case~\citep[Proposition
6.2.8]{MT93}).
\end{proof}

Under A\ref{Ex1:A1}-\ref{Ex1:A2}, it si easily checked that any twice
continuously differentiable function $V : \Rset^n \to \Rset$ is in the domain
of $\gen$ and $LV(x) = \gen V(x) $ for all $x \in \Rset^n$. Observe indeed that
$t \mapsto LV(X_t)$ is integrable $\pr_x$-a.s. for all $x \in \Rset^n$ and $t
\mapsto V(X_t) - V(X_0) - \int_0^t LV(X_s) ds$ is a right-continuous local
martingale. Hence $V \in \mathcal{D}(\gen)$ and $LV= \gen V$. We now establish
drift inequalities for different test functions $V$.
\begin{prop}
  \label{prop:Ex1_DriftIneq}
  Assume A\ref{Ex1:A1}-\ref{Ex1:A2}.  Let $V: \Rset^n \to [1,+\infty)$ be a
  twice continuously differentiable function such that $ V(x) = \exp(\iota \ 
  |x|^m)$ outside a compact set, for some $0<m<1$ and $\iota >0$. Then
  $\sup_{|x| \leq M} \gen V(x) < \infty$ and for all $|x| \geq M$,
  \begin{enumerate}[(i)]
  \item If $0<m<1-p$,
\[ 
\gen V(x) \leq - \iota^{\frac{1+p}{m}} m r \left[ \ln V(x)\right]^{1 - (\frac{1+p}{m})} V(x) \ \ \left(1 +  o(1) \right);
\]
\item If $m=1-p$,
\[
\gen V(x) \leq - \iota^{\frac{1+p}{1-p}} (1-p) \left\{r - (1/2) \lambda_+ \iota (1-p) \right\} \   \left[ \ln V(x)\right]^{ - 2\frac{p}{1-p}} V(x) \ \ \left(1 +  o(1) \right).
\]
  \end{enumerate}
\end{prop}
\begin{proof}
  Under the stated assumptions, $\sup_{ \{x, |x| \leq M \}} \gen V(x) < \infty$.
  By definition of $\gen$, we have for all $|x| \geq M$,
\[
 \gen V(x) \leq -\iota m \left( r - (1/2) \lambda_+ \iota m |x|^{p+m-1}\right) |x|^{m-1-p} V(x) + (1/2) \iota m n \Lambda  |x|^{m-2} V(x).
\]
\end{proof}
As a direct application of Proposition~\ref{prop:CS_PiInteg} and Theorem
\ref{theo:ContMax}(\ref{ContMax_Csq2}), we have 
\begin{theo}
  \label{theo:Ex1_OurResults1}
Assume A\ref{Ex1:A1}-\ref{Ex1:A2}.
\begin{enumerate}[(i)]
\item For all $\iota>0$ such that $r - (1/2) \lambda_+ \iota (1-p)>0$, 
\[
\int  \pi(dx)  \exp(\iota |x|^{1-p}) < \infty,
\]
where $\pi$ is the invariant probability distribution of the Markov process that solves (\ref{eq:StochEquDiff}).
\item  There exists a closed petite set $C$ such that for any $0<m < 1-p$, $0<\iota_1<\iota_2$ and $\delta>0$, there exists a finite constant $c$ such that 
  \begin{equation}
\label{eq:Ex1_MomExpo}
 \esp_x \left[ \exp(\iota_1 \;  \{ \tau_C(\delta)\}^{\frac{m}{1+p}}) \right] \leq c \exp(\iota_2 |x|^m).
  \end{equation}
If $m=1-p$, (\ref{eq:Ex1_MomExpo}) still holds for any $0<\iota_1<\iota_2$ such
that $r - (1/2) \iota_2 \lambda_+  (1-p)>0$.
\end{enumerate} 
\end{theo}

  The results of Theorem~\ref{theo:Ex1_OurResults1} can be compared to those
  by~\citep{Malyshkin:2001}, where subexponential ergodicity in total variation norm of a
  diffusion satisfying the conditions A\ref{Ex1:A1}-\ref{Ex1:A2} is addressed.
  The technique used in~\citep{Malyshkin:2001} is based on the coupling method.
  Theorem~\ref{theo:Ex1_OurResults1}(i) states the same result as~\citep[Lemma
  3]{Malyshkin:2001}. Nevertheless, Theorem~\ref{theo:Ex1_OurResults1}(ii) yields a
  stronger control of delayed return-time to a closed petite set than those
  obtained in~\citep[Theorem 5]{Malyshkin:2001}. They show that for all $0< \alpha <(1/2)
  (1-p)$ there exists a constant $c_\alpha$ such that
\[
\esp_x \left[ \exp( \tau_C(\delta)^{\alpha}) \right] \leq c_\alpha \exp(|x|^{2 \alpha}),
\]
and this remains valid for $\alpha=(1-p)/2$ if $r - (1/2) \lambda_+ (1-p)>0$.
Theorem~\ref{theo:Ex1_OurResults1}(ii) claims that for all $0 < \alpha < (1-p)
(1+p)^{-1}$ and $\iota>1$, $\esp_x \left[ \exp( \tau_C(\delta)^{\alpha})
\right] \leq c_\alpha \exp(\iota \; |x|^{(1+p)\alpha})$
and for $\alpha= (1-p) (1+p)^{-1}$,  $\esp_x \left[ \exp( \iota_1 \tau_C(\delta)^{\alpha}) \right] \leq c_\alpha \exp(\iota_2 \; |x|^{(1+p)\alpha}) $ for all $0<\iota_1<\iota_2$ such that $r - (1/2) \; \iota_2 \lambda_+ (1-p)>0$. \\
As a direct application of Theorem~\ref{theo:FR_ergod}, we obtain the following
results for $f$-ergodicity at a subgeometric rate.
\begin{theo}
  \label{theo:Ex1_OurResults2}
  Assume A\ref{Ex1:A1}-\ref{Ex1:A2} and let $\pi$ be the invariant probability
  distribution of the Markov process that solves (\ref{eq:StochEquDiff}).  Then the process is subgeometrically $f$-ergodic: 
for
  any $x \in \Rset^n$, the limits~(\ref{eq:Ergodicite}) to (\ref{eq:ControleExplicite2}) hold with
$V(x) \sim \exp( \iota |x|^{1-p})$ for some positive $\iota$ such that $r - 0.5
  \lambda_+ \iota (1-p)>0$, $f_\ast(x) \sim |x|^{-2p} \exp(\iota |x|^{1-p})$
  and $r_\ast (t) \sim t^{-2p/(1+p)} \exp(\{\iota' t\}^{(1-p)/(1+p)})$ where
\[
\iota' = \iota^{\frac{1+p}{1-p}} (1+p) \left\{r - (1/2) \lambda_+ \iota (1-p)
\right\}.
\]
  \end{theo}
  In~\citep{Malyshkin:2001}, only the  convergence in total variation norm of the semi-group  $\{P^t \}_{t \geq 0}$ to the invariant probability $\pi$ is addressed: is is established  that the process is ergodic at the rate $r_\ast^{\mathrm{M}}(t) \propto \exp(\delta t^{(1-p)/2})$ for some $\delta >0$, and in that case, the dependence upon the initial point in (\ref{eq:Ergodicite}) is $V^{\mathrm{M}}(x) \sim \exp(\delta |x|^{1-p})$. Theorem~\ref{theo:Ex1_OurResults2}  improves these results and  also provides rates of convergence in $f$-norm for unbounded functions $f$. \\
  We reported in Theorem~\ref{theo:Ex1_OurResults2} the values $(V,f_\ast,
  r_\ast)$ that yield the best rate of convergence in total variation norm.
  Proposition~\ref{prop:Ex1_DriftIneq} shows that one could establish the drift
  inequality~(\ref{eq:gendrift}) with $V(x) \sim \exp( \iota |x|^m)$ for some
  $0< m <1-p$; this would imply the limits~(\ref{eq:Ergodicite}) to
  (\ref{eq:ControleExplicite2}) with $V(x) \sim \exp(\iota |x|^m)$, $f_\ast(x)
  \sim |x|^{m-1-p} \exp(\iota |x|^m)$ and $r_\ast(t) \sim t^{(m-1-p)/(1+p)}
  \exp(\iota' |x|^{m/(1+p)})$ for all $0< \iota' < \iota$. We thus obtain a
  weaker maximal rate function $r_\ast$, and a weaker maximal norm $\| \cdot
  \|_{f_\ast}$, but this has to be balanced with the fact that the dependence
  upon the initial value (\ie the quantity $V(x)$) is weaker too. Similarly,
  polynomially increasing controls $V(x)$ could be considered, thus limiting
  the rate $r_\ast$ (resp. the function $f_\ast$) to the class of the
  polynomially increasing rate functions (resp. to the class of the
  polynomially increasing function). These discussions illustrate the fact that
  the pair $(\phi,V)$ that solves (\ref{eq:gendrift}) is not unique, and this
  results in balancing the pair $(r_\ast, f_\ast)$ and the dependence upon the
  initial value $x$.

%%%%%%%%%%%%%%%%%%%
%% langevin
%%%%%%%%%%%%%%%%%%%

\subsection{Langevin tempered diffusions on $\Rset^n$}
\label{sec:Ex2}
Let $\pi : \Rset^n \to  (0,\infty)$ satisfying
\debutB
\item \label{Ex2:B1} $\pi$ is, up to a normalizing constant,  a positive and thrice continuously differentiable density on $\Rset^n$, with respect to the Lebesgue measure.
\finB
Let $\sigma(x)=|\ln \pi(x)|^d$ for some $d>0$ and define the diffusion matrix by $a(x) = \sigma^2(x) \Id_n$, and the drift vector by  $b(x) =(b_1(x) ,\cdots, b_n(x))'$ where
\[
b_i(x) = (1/2) \sum_{j=1}^n a_{ij}(x) \ \partial_{x_j} \log \pi(x) + (1/2)  \sum_{j=1}^n \partial_{x_j} a_{ij}(x), \qquad 1 \leq i \leq n.
\]
Observe that since $\pi$ is defined up to a normalizing constant, we can assume that $\sigma(x)>0$ for all $x$.
Our objective is to study the ergodicity of the solution to the   stochastic integral equation
\begin{equation}
  \label{eq:Ex2_StochEqu}
  X_t= X_0 + \int_0^t b(X_s) ds + \int_0^t \sigma(X_s) dB_s
\end{equation}
where  $\{B_t\}_t$ is a $n$-dimensional Brownian motion.  This diffusion is the so-called Langevin diffusion and the drift vector $b$ is defined in such a way that $\pi$ is, up to a multiplicative constant, the density of the unique invariant probability distribution. Note that this model is not a particular case of the elliptic diffusion of section \ref{sec:Ex1} since here, $\sigma$ may be an unbounded function ($\sigma =|\ln \pi(x)|^d$).   \\
Fort and Roberts investigate the behavior of these diffusions when $\pi$ is
polynomially decreasing in the tails and address ergodicity in total variation
norm and in $f$-norm as well~\citep{Fort:Roberts:2005}. They consider the case $\sigma(x) =
\pi^{-d}(x)$ ($d>0$) and show that the rate of convergence in total variation
norm and in $f$-norm for $ f(x) \sim \pi^{-\kappa}(x)$ ($\kappa>0$) depends on
$d$.  When $d$ is lower than some critical temperature $d_\ast$, the process is
ergodic at a polynomial rate, and when $d$ is larger than $d_\ast$, the process
is uniformly ergodic in total variation norm and geometrically ergodic
otherwise~\citep[Theorem 16]{Fort:Roberts:2005}. Fort and Roberts thus proved that the rate
of convergence can be improved by choosing a diffusion coefficient $\sigma$
which is small when the process is close to the modes of $\pi$ and big when it
is far from the modes.  The objective of this section is to investigate the
case when $\pi$ is super-exponentially decreasing in the tails.  We assume that
\debutB
\item \label{Ex2:B2} there exists  $0 < \beta<1$ such that for all
 large $|x|$,
  \begin{eqnarray*}
     && |x|^{1-\beta} \;  \pscal{\partial \ln \pi(x)}{\frac{x}{|x|}}< 0, \\
    && 0 < \liminf_{x \to \infty} |\partial \ln \pi(x)|| \ln \pi(x)|^{1/\beta-1}
 \leq  \limsup_{x \to \infty} |\partial \ln \pi(x)| | \ln \pi(x)|^{1/\beta-1} < \infty,  \\
  && \limsup_{x \to \infty} \mathrm{Tr}\left(\partial^2 \ln \pi(x)\right)|\partial \ln \pi(x)|^{-2} =0.
  \end{eqnarray*}
  \finB 
 The class of
  density $\pi$ described by B\ref{Ex2:B1}-\ref{Ex2:B2} contains densities that
  are super-exponential in the tails. The Weibull distribution on $(0,\infty)$
  with density $\pi(x) \propto x^{\beta-1} \exp(-\alpha x^\beta)$ satisfies
  B\ref{Ex2:B2}. For multidimensional examples, see
  e.g.~\citep{Kluppelberg:1988,Roberts:Tweedie:1996,Fort:Moulines:2000}.  Following the same steps as in
  Section~\ref{sec:Ex1}, we can prove that under B\ref{Ex2:B1}-\ref{Ex2:B2} and
  provided the process is regular, there exists a solution to
  (\ref{eq:Ex2_StochEqu}) which is an almost surely continuous stochastic
  process and is unique up to equivalence. This solution is an homogeneous
  strong Markov process whose transition functions are Feller functions. Under
  B\ref{Ex2:B2}, the process is regular whatever $d>0$; this can
  be proved as in the previous section (by choosing $V = 1 + \pi^{-2}$, see (\ref{eq:Ex2:DriftV}) below). \\
  These assumptions also imply that $\pi$ is (up to a scaling factor) the
  density of an invariant distribution of the diffusion process, any skeleton
  chain is $\psi$-irreducible and
  compact sets are closed  petite sets   (\citep[Proposition 15]{Fort:Roberts:2005}).  \\
  Let $V: \Rset^n \to [1,\infty)$ be a twice-continuously differentiable
  function such that $V(x) = 1 + \pi^{-\kappa}(x)$ outside a compact set; then
  $ \gen V(x) = LV(x) = \ell_1(x) + \ell_2(x)$ where $L$ is the diffusion
  operator (\ref{eq:Ex1:Loperator}) and for large $|x|$,
\begin{equation}
  \label{eq:Ex2:DriftV}
  \ell_1(x) =  - \frac{\kappa (1-\kappa)}{2} \; \frac{\pi^{-\kappa}(x)}{1+
 \pi^{-\kappa}(x)} \; \left(\frac{ |
\partial \ln \pi(x)|  }{|\ln \pi(x)|^{1-1/\beta}}  \right)^2   \  \ |\ln \pi(x)|^{2(d+1-1/\beta)} \ V(x), 
\end{equation} 
and $\ell_2(x) = o(\ell_1(x))$. In~\citep[Theorem 3.1]{Stramer:Tweedie:1999}, it is established that 
the process is geometrically ergodic if and only if  $d \geq  1/\beta-1$. From
(\ref{eq:Ex2:DriftV}), we are able to retrieve these results and we also  prove that
when $ 0 \leq d <   1/\beta-1$, the process is subgeometrically ergodic.
Observe indeed that for large $|x|$, (\ref{eq:Ex2:DriftV}) and B\ref{Ex2:B2}
imply 
\[
\gen V(x) \leq - c_\kappa \left[ \ln V(x) \right]^{-\alpha} V(x), \quad \text{where}
\; \alpha=2(1/\beta-1-d) , \quad \text{and} \;  c_\kappa> 0 \Longleftrightarrow 0< \kappa <1. 
\]
Hence, if $\alpha \leq 0$, the process is $V$-geometrically
ergodic~\citep[Theorem 6.1]{Meyn:Tweedie:1993} (see also section~\ref{sec:Generator}); if $\alpha >0$, it is subgeometrically
ergodic as a consequence of Theorems \ref{theo:FR_ergod} and \ref{theo:generator2df}. \\
A polynomially increasing drift function can also be considered: we can assume without loss of generality that for large $x$, $\ln \pi(x) < 0$  since $\pi$ is defined up to a multiplicative constant. We thus set $V(x) = 2 +\mathrm{sign}(\kappa) \ ( - \ln \pi(x) )^\kappa$ outside a compact set. Then  for large $x$, 
\[
\gen V(x) \leq -\frac{ |\kappa|}{2}  \left( - \ln \pi(x) \right)^{\kappa+1+2(d-1/\beta)} \ \left(\frac{ |
\partial \ln \pi(x)|  }{|\ln \pi(x)|^{1-1/\beta}}  \right)^2 \ \ \left(1 + o(1) \right),
\]
and there exists a constant $c>0$ such that for large $x$, 
\begin{equation}
  \label{eq:DriftEx2}
   \gen V(x) \leq - c V^{1- \alpha}(x), \quad \text{where}\; \alpha=2\kappa^{-1}(1/\beta-d-(1/2)).
\end{equation}
First consider the case when $\kappa>0$. If $1/\beta -1 < d <1/\beta -(1/2)$,
the drift condition (\ref{eq:DriftEx2}) and Theorems~\ref{theo:FR_ergod} and
\ref{theo:generator2df} yield polynomial ergodicity. For example, this implies
convergence in total variation norm at the rate $r(t) \sim t^{1/\alpha-1}$. If
$d=1/\beta -(1/2)$, then $\alpha=0$ and the process is geometrically ergodic.
In the case when $\kappa$ can be set negative and $1-\alpha>0$ \ie\ when
$d>1/\beta -(1/2)$, the process is uniformly ergodic: there exist $\lambda<1$
and a constant $c<\infty$ such that for all $x$,
\[
\lambda^{-t} \ \| P^t(x,\cdot) - \pi(\cdot) \|_{\mathrm{TV}} \leq c,
\]
and the convergence does not depend on the starting point. \\
The above discussions are summarized  in  the following
theorem. The first  part (resp. third part) results from \citep[Theorem 3.1]{Stramer:Tweedie:1999}
(resp. \citep[Theorem 6.1]{Meyn:Tweedie:1993}). The second assertion is a consequence of
Theorem~\ref{theo:FR_ergod}. The last assertion was already proved by~\citep[Theorem 3.1]{Stramer:Tweedie:1999} for one-dimensional diffusions ($n=1$).
\begin{theo}
  \label{theo:Ex2}
  Consider the Langevin diffusion on $\Rset^n$ solution to the equation
  (\ref{eq:Ex2_StochEqu}) where the target distribution $\pi$ satisfies
  B\ref{Ex2:B1}-\ref{Ex2:B2}.
  \begin{enumerate}[(i)]
  \item If $0 \leq d < 1/\beta -1$, the process fails to be geometrically ergodic.
\item If $0 \leq d <  1/\beta -1$, the process is subgeometrically ergodic:  the limits (\ref{eq:Ergodicite}) to (\ref{eq:ControleExplicite2}) hold with $V(x) \sim \pi^{-\kappa}(x)$, $f_\ast(t)  \sim \pi^{-\kappa}(x) \left| \ln \pi(x) \right|^{-2(1/\beta-1-d)}$ and $\ln r_\ast(t) \sim c_\kappa t^{\beta/(2-\beta-2 d \beta)}$ for all $0 < \kappa <1$.
\item If $d \geq 1/\beta -1$, then for all $0 < \kappa <1$, the diffusion is
  $V$-geometrically ergodic with $V(x) = 1+ \pi^{-\kappa}(x)$.
\item If $d >  1/\beta -(1/2)$, the diffusion  is uniformly ergodic.
  \end{enumerate}
\end{theo}
This theorem extends earlier results to the multi-dimensional case and provides subgeometrical rates of convergence of the 'cold' Langevin diffusion, for a wide family of norms. 
We established that for a given $\pi^{-\kappa}$-norm, the minimal rate of convergence is achieved with $d=0$ and in that case, the rate coincides with the rate of convergence of the symmetric random-walk Hastings-Metropolis algorithm (\citep[Theorem 3.1]{Douc:Fort:Moulines:Soulier:2004}). This rate can be improved by choosing a diffusion matrix which is heavy where $\pi$ is light and conversely. When $d$ is larger than the critical value $d_\ast = 1/\beta-1$, the process is geometrically ergodic; when $d$ is lower that $d_\ast$, the process can not be geometrically ergodic and we prove that it is subgeometrically ergodic. The conclusions of Theorem~\ref{theo:Ex2} are similar to those of ~\citep[Theorem 16]{Fort:Roberts:2005}, that address the case when $\pi$ is polynomial in the tails.  \\
We assumed that $\sigma = | \ln\pi|^d$.  A first extension is to consider a sufficiently smooth function $\sigma$ such that $\sigma(x) \sim | \ln \pi(x)|^d$ for large $|x|$; this yields similar conclusions and details are omitted. A second extension consists in  the case when $\sigma(x) \sim \pi^{-d}(x)$. In this latter case, following the same lines, it is easily verified that for small enough values of $d$, the process is regular (the set of the admissible values is in the range $(0,1/2]$), and the process is $V$-geometrically ergodic with a test function $V(x) \sim \pi^{-\kappa}(x)$, $\kappa>0$. Details are omitted and left to the interested reader.

%%%%%%%%%%%%%%%
%% hypoellitpic
%%%%%%%%%%%%%%%

\subsection{Stochastic damping Hamiltonian system}

Both examples of the previous sections assumed that the diffusion process is
elliptic. However the drift condition (\ref{eq:gendrift}) enables us to
consider also hypoelliptic diffusion that we will illustrate on the example of
a simple stochastic damping Hamiltonian system, i.e. let $x_t$ (resp. $y_t$) be
the position (resp. the velocity) at time $t$ of a physical system moving in
$\Rset^n$
\begin{equation}\label{eq:ham}
\begin{array}{l}
dX_t=Y_tdt\\
dY_t=\Sigma(X_t,Y_t)dB_t-(c(X_t,Y_t)Y_t+\partial_x U(X_t))dt
\end{array}
\end{equation}
where $-\partial_x U$ is some friction force, $-c(x,y)y$ is the damping force and $\Sigma(x,y)dB$ is  a random force where $(B_t)$ is a standard Brownian motion in $\Rset^n$. This system has been studied from the large and moderate deviations point of view by Wu \citep{Wu:2001} where he also establishes the exponential ergodicity under various set of assumptions.

As our goal is not to consider the model in its full generality but to illustrate the subexponential behavior of hypoelliptic diffusion, via the simple use of drift condition (\ref{eq:gendrift}), hereafter we will consider the particular (but also current in practice) case where the damping and random forces are constant $c(x,y)=c \ \Id_n$ and $\Sigma(x,y)=\sigma \ \Id_n$, $c$ and $\sigma$ being positive constants (as, if one is identically equal to 0, there is none of the usual ergodic properties such as positive recurrence). We will assume moreover that the potential $U$ is lower bounded and continuously differentiable over $\Rset^n$. In this case, the system is known to have an unique invariant measure given by
$$\pi(dx,dy)=e^{-{2c\over\sigma}H(x,y)}dxdy$$
where $H$ is the Hamiltonian
given by $H(x,y)={1\over 2}|y|^2+U(x).$

Let us first ensure the existence of solutions and aperiodicity for the process
$Z_t=(X_t,Y_t)$ via the following proposition due to Wu \citep[Lemma 1.1,
Proposition 1.2]{Wu:2001}
\begin{prop}
  For every initial state $z=(x,y)\in \Rset^{2n}$, the SDE (\ref{eq:ham})
  admits an unique weak solution $\mathbf{P}_z$ which is non explosive.
  Moreover denoting $(P^t(z,dz'))_t$ the associated semi group of transition,
  we have that for every $t>0$ and every $z\in\Rset^{2n}$,
  $P^t(z,dz')=p_t(z,z')dz'$ and $p_t(z,z')>0,$ $dz'-a.e.$ The density
  $p_t(z,\cdot)$ is moreover continuous, and the process is thus strongly
  Feller.
\end{prop}
As a consequence, the solution is a strong Markov process,  all the skeletons are irreducible and compact sets are petite sets.\\
Let us build an example of polynomially ergodic stochastic damping Hamiltonian
system in dimension $1$. We rewrite the system as
\begin{equation}\label{eq:ham2}
\begin{array}{l}
dX_t=Y_tdt\\
dY_t=\sigma dB_t-(c Y_t+ U'(X_t))dt,
\end{array}
\end{equation}
and assume that $U$ is $C^2$, and there exist $0<p<1$
and positive constants $a,b$ such that for $|x|$ large enough
\begin{equation}\label{eq:condhamil}
a|x|^{p-1}\le U'(x)\le b |x|^{p-1}.
 \end{equation}
 The fact that $p$ is less than 1 implies that $(Z_t)_{t \geq 0}$ cannot be
 exponentially ergodic~\citep[Theorem 5.1]{Wu:2001}.  We now exhibit a drift
 function satisfying (\ref{eq:gendrift}). Consider positive constants
 $\alpha,\beta$ and a smooth positive function $G$ such that for $m$,
 $1-p<m\le1$, $G'(x)=|x|^m$ for large $|x|$; define a twice continuously
 differentiable function $V\geq 1$ such that for large $x,y$,
$$V(x,y)=\alpha (y^2/2+U(x))+\beta ( G'(x) y+ c G(x)).$$
% and note $V=v-\inf v +1.$
By definition of $\gen$, it holds
$${\gen}V_m(x,y)={1\over 2}\sigma^2 \ \partial^2_yV_m(x,y) +y \  \partial_x V_m(x,y) -(cy+U'(x)) \ \partial_y V_m(x,y)$$
so that
\begin{eqnarray*}
{\gen}V_m(x,y)&=&{1\over 2}\alpha\sigma^2+y(\alpha U'(x)+\beta G''(x)y+\beta c G'(x))-(cy+U'(x))(\alpha y+\beta G'(x))\\
&=&{1\over 2}\alpha\sigma^2+(\beta G''(x)-c\alpha)y^2-\beta G'(x)U'(x).
\end{eqnarray*}
Fix $\delta <0$; since $m \leq 1$, we choose $\beta$ small enough so that
$\beta G''(x)-c\alpha<\delta<0$ for all large $x$. Furthermore, for all large
$|x|$, $G'(x)U'(x)\ge b|x|^{p-1+m}$. Hence, there exist positive constants
$K,L$ such that
$${\gen }V_m(x,y)\le K - L \; V_m(x,y)^{p-1+m\over m+1}.$$
Condition
(\ref{eq:gendrift}) holds with $\phi_m(v)\propto v^{p-1+m\over m+1}$ and
${p-1+m\over m+1}<1$. Application of the results of Section~\ref{sec:ergodic}
now implies that the process $(Z_t)_{t \geq 0}$ is polynomially-ergodic.
\par

Let $k\ge1$ and define a twice continuously differentiable function $V_{m,k} \geq 1$ such
that for large $x,y$
$$V_{m,k}(x,y)=V_m^k(x,y).$$
Then for large $x,y$, the above calculations  yield
\begin{eqnarray*}
\gen V_{m,k}(x,y)&=&(\gen V_m(x,y))V_m^{k-1}(x,y)+{1\over 2}\sigma^2 (\partial_y V_m(x,y))^2V_m^{k-2}(x,y)\\
&=&\left( \gen V_m(x,y))+{1\over 2}\sigma^2 {(\partial_y V_m(x,y))^2\over V_m(x,y)}\right) V_m^{k-1}(x,y)\\
&\le & (K' - L V_m(x,y)^{p-1+m\over m+1})V_m^{k-1}(x,y)\\
&\le & K''-L' V_m^{{p-1+m\over m+1}+k-1}
\end{eqnarray*}
for some positive constant $K', K'', L'$. This inequality is once again the
condition (\ref{eq:gendrift}) with $\phi_{m,k}(v)=v^{({p-2\over
    m+1}+k)k^{-1}}$. These discussions are summarized in the following
Theorem.
\begin{theo}
\label{theo:hypo}
Let $U$ be a twice continuously differentiable function, lower bounded on
$\Rset$ satisfying (\ref{eq:condhamil}) for some $0<p<1$. Then $(Z_t)_{t \geq
  0}$ is not exponentially ergodic but is polynomially ergodic~: for any $m$
such that $1-p< m\le 1$ and any $k \geq 1$, the limits
(\ref{eq:Ergodicite}-\ref{eq:ControleExplicite2}) hold with $V_{mk}$ defined
above, $\phi_{m,k}(v) \propto v^{({p-2\over m+1}+k)k^{-1}}$, $f_*=\phi_{m,k}\circ V_{m,k}$
and $r_*(t) \propto t^{{k(m+1)\over 2-p}-1} $.
\end{theo}
Observe that the process $(Z_t)_{t \geq 0}$ is polynomially ergodic at any
order and we strongly believe it is subexponentially ergodic. This sub
exponential case is left to the interested reader. The multidimensional case
is more intricate in the
choice of the drift function and we do not pursue here in this direction. \\
This example shows that our conditions are sufficiently flexible to consider
the hypoelliptic diffusions as well as the elliptic ones.

%%%%%%%%%%%%%%%%
%% eds saut
%%%%%%%%%%%%%%%%

\subsection{Compound Poisson-process driven Ornstein-Uhlenbeck process}
\label{sec:Ex3}
In this section we consider an example of Fort-Roberts \citep{Fort:Roberts:2005} where subgeometric ergodicity can be achieved where they only obtain polynomial ergodicity. Let us first recall the model. Let $X$ be an Ornstein-Uhlenbeck process driven by a finite rate subordinator:
$$dX_t=-\mu X_t+dZ_t$$
and $Z_t=\sum_{i=1}^{N_t}U_i$, where $(U_i)_{i\ge1}$ is a sequence of i.i.d.r.v. with probability measure $F$, and $(N_t)$ is an independent Poisson process of rate $\lambda$. We suppose the recall coefficient $\mu$ to be positive. Remarking that only when $F$ is sufficiently (even extremely) heavy tailed, $X$ fails to be exponentially ergodic, Fort-Roberts \citep{Fort:Roberts:2005} give conditions for which $X$ is polynomially ergodic. Namely, denote $G$ the law of the log jump sizes ($G(A)=F(e^A)$), and assume that for all $\kappa>0$, $\int e^{\kappa x}dG(x) = +\infty$. Lemma 17 of Fort-Roberts then prove that $X$ is not exponentially ergodic and give examples where $X$ is positive recurrent and polynomially ergodic, namely when for some $r>1$, $\int_0^\infty [\log(1+u)]^r F(du)$ is finite. Such assertion may be useful considering
$$F(dx)= {C_k^{-1}\over x(\log(x))^k}dx\qquad \qquad k>1$$
$$F(dx)={C_{\beta,c}^{-1}e^{-c(\log(x))^\beta}\over x}dx\qquad\qquad \beta\le1.$$

We shall strengthen their result by

\begin{prop}
Suppose that $(X_t)$ is aperiodic and that for some $\delta<1$, $\alpha>0$
$$\int_0^\infty e^{\alpha (\log(1+x))^\delta} F(dx)<\infty.$$
Then, the conclusions of Theorem \ref{theo:FR_ergod} hold with $V(x)=e^{\alpha'(\log x)^{\delta'}}$ (and $\alpha'<\alpha$ if $\delta'=\delta$), and $\phi(v)=v^{(1-\delta')/\delta'}$,  $r_\ast(t)=\mathfrak{a}t^{-(1+\delta')}e^{\mathfrak{b}t^{\delta'}}$, $f_\ast=\phi\circ V$.
\end{prop}

{\it Proof}. We shall use the drift conditions introduced previously for the generator defined by for all functions $V$ in the extended domain of the generator
$${\cal A}V(x) =\lambda \int_0^\infty (V(x+u)-V(x)) F(du)-\mu xV'(x).$$
Choosing $V(x)=(\log(x))^r$, as in Fort-Roberts \citep[Lemma 18]{Fort:Roberts:2005}, for sufficiently large $x$ ensures the polynomial ergodicity at the previous rate. Consider now $V(x)=e^{\alpha'(\log x)^{\delta'}}$, so that
\begin{eqnarray*}
{\cal A}V(x)&=&\lambda  \int_0^\infty (e^{\alpha'(\log (x+u))^{\delta'}}-e^{\alpha'(\log x)^{\delta'}})F(du)-\alpha'\delta'\mu {e^{\alpha'(\log x)^{\delta'}}\over (\log x)^{1-\delta'}}\\
&\le& -\alpha'^{2-\delta'}\mu\delta'{V\over (\log V)^{(1-\delta')/\delta'}} +b
\end{eqnarray*}

recalling that % all bounded sets are petite in this example and that
 for large $x$
$$e^{\alpha'(\log (x+u))^{\delta'}}-e^{\alpha'(\log x)^{\delta'}}\sim \delta'    {e^{\alpha'(\log x)^{\delta'}}\over (\log x)^{1-\delta'}}\log(1+u/x);$$
the dominated convergence theorem ends the argument.

\section{Skeleton chain and moderate deviations}\label{discret}

We consider here an important field of application for this subgeometric rate,
namely moderate deviations for bounded additive functionals of Markov process.
In fact, Proposition \ref{prop:lienContDiscret} gives us more than a way to
deal with subexponential ergodicity, it also implies a drift condition in the
sense of Douc-Fort-Moulines-Soulier \citep{Douc:Fort:Moulines:Soulier:2004}
which will enables us, at least in a bounded test function framework, to extend
to the continuous time case some limit theorems tailored for the subexponential
regime by Douc-Guillin-Moulines \citep{Douc:Moulines:Guillin:2005} such as moderate deviations.
Moderate deviations are concerned with the asymptotic for centered $g$ with
respect to $\pi$ and for $0\le t\le T$ of
$$
S^\epsilon_t={1\over \sqrt{\epsilon}h(\epsilon)}\int_0^tg(X_{s/\epsilon})ds
$$
where as $\epsilon$ tends to 0, $h(\epsilon)\to\infty$ but $\sqrt{\epsilon}h(\epsilon)\to0$, namely a regime between the large deviations and the central limit theorem. We may then state (proofs will be done in appendix.
)

\begin{theo}\label{thm:mdp}
  Assume that $\mathbf{D(C,V,\phi,b)}$ holds with $\sup_C V< \infty$, and some
  skeleton chain is $\psi$-irreducible.
\begin{enumerate}[(i)]
\item   For all $m>0$, there exist a function $ W : \sfX \to [\phi(1), \infty)$, a small set $\tilde C$ for the skeleton $P^m$ and a positive constant $b'$  such that $\sup_{\tilde C}W$ is finite, and on $\sfX$,
$$
P^m W \leq W -\phi \circ W + b' \ind_{\tilde C}, \qquad  \text{and} \qquad \phi \circ V  \leq W \leq \kappa V.
$$
\item Assume that $X_0$ is distributed as $\mu$ and $\mu(V)<\infty$ and that $g$ is a bounded mapping from ${\sfX}$ to $\Rset^n$ with $\pi(g)=0$. Suppose moreover that for all positive $a$
$$\lim_{\epsilon\to 0}{1\over h^2(\epsilon)}\log\left( \epsilon H_\phi^{-1}\left({a~h(\epsilon)\over\sqrt{\epsilon}}\right)\right)=\infty$$

then $\pr_\mu\left(S^\epsilon_\cdot\in \cdot\right)$ satisfies a moderate deviation principle in $C_0([0,1],\Rset^n)$ (the space of continuous functions from $[0,1]$ to $\Rset^n$ starting from 0) equipped with the supremum norm topology, with speed ${1\over h^2(\epsilon)}$ and rate function $I_g^h$, i.e. for all Borel set $A\in C_0([0,1],\Rset^n)$
\begin{eqnarray*}
{\rm -}\inf_{\gamma \in int(A)} I_g^h(\gamma)&\le &\liminf_{\epsilon\to0}{1\over h^2(\epsilon)}\log\pr_\mu\left(S^\epsilon_\cdot\in A\right)\\
&\le& \limsup_{\epsilon\to0}{1\over h^2(\epsilon)}\log\pr_\mu\left(S^\epsilon_\cdot\in A\right)\le{\rm -}\inf_{\gamma \in cl(A)}I_g^h(\gamma)
\end{eqnarray*}
where $I_g^h$ is given by
\begin{equation}
I_g^h(\gamma):=\left\{\begin{array}{ll}
 \displaystyle{1\over 2}\int_0^1\sup_{\zeta\in\Rset^n}\left\{\langle \dot\gamma(t),\zeta\rangle-{1\over 2}\sigma^2(\langle g,\zeta\rangle)\right\}dt&\mbox{\rm if }d\gamma(t)\mbox{\rm =}\dot\gamma(t)dt,~\gamma(0)\mbox{\rm =}0,\\ \\
+\infty& \mbox{\rm else},\end{array}\right.
\end{equation}
and 
\begin{equation}
\begin{array}{ll}
\sigma^2(\langle g,\zeta\rangle)&\displaystyle=\lim_{n\to\infty}{1\over n}\esp_\pi\left(\int_0^n g(X_s)ds\right)^2\\ \\
&\displaystyle=2\int_\sfX\langle g,\zeta \rangle\int_0^\infty P_t \langle g,\zeta \rangle dt~ d\pi.\end{array}
\end{equation}
\end{enumerate}
\end{theo}

The proof is in Section \ref{proof:mdp}. 

To the authors' knowledge, this moderate deviations result (even for bounded
function) is the first one for Markov processes which are not exponentially
ergodic. It extends then results of Guillin \citep[Th 1.]{Guillin:2001} or Wu \citep[Th.
2.7]{Wu:2001} in the subexponential setting. As expected, all ranges of speed
are not allowed for such a theorem but are limited by the ergodicity of the
process (we refer to Douc-Guillin-Moulines \citep[Sect. 4]{Douc:Moulines:Guillin:2005} for a complete
discussion on this interplay). The extension of this moderate deviation
principle to unbounded function is left for further research, as well as
extension to inhomogeneous functional and averaging principle, those subjects
needing particular tools and developments.
 
\appendix

\section{Proofs}
\label{app:Proofs}
\subsection{Proof of Theorem~\ref{theo:ContMax}}
\begin{lem} \label{lem:cqdlag}
  For any $M >0$ and for any cad-lag function $g$,
\begin{equation}\label{eq:cadlag}
\lim_{\epsilon \to 0}\sum_{k=1}^{\lfloor M/\epsilon \rfloor} \left|\int_{\epsilon (k-1)}^{\epsilon k} (g(s)-g(t_{k-1}))ds \right| = 0.  
\end{equation} 
\end{lem}
\begin{proof}
First note that $g$ is bounded since it is a cad-lag function. Let $\eta>0$ be an arbitrary real. For any $x \in [0,M]$, there exists an interval $(x-\alpha,x+\alpha)$ such that 
$$
\forall s \in (x-\alpha,x), \ |g(s)-g(x -)| < \eta/2 \quad \mbox{and} \quad \forall s\in [x,x+\alpha), |g(s)-g(x)| < \eta/2
$$ 
Thus, for any $(u,v)$ in $(x-\alpha,x)\times (x-\alpha,x) $ or in $ [x,x+\alpha) \times  [x,x+\alpha)$, $|g(u)-g(v)|\leq \eta$.  
By compacity of $[0,M]$, there exists a finite number  $M_\eta$ of such intervals $(x_i-\alpha_i, x_i+\alpha_i)$ which covers $[0,M]$.  Taking $\epsilon$ sufficiently small, it can be easily checked that any interval $[\epsilon (k-1), \epsilon k]$ is included in some interval $(x_i-\alpha,x_i+\alpha)$. Now, if some $x_i \in [\epsilon (k-1), \epsilon k]$, write $\sup_{u,v \in [\epsilon (k-1), \epsilon k]}|g(u)-g(v)|\leq 2 \sup_{x \in [0,M]}|g(x)|$. Otherwise, we have $\sup_{u,v \in [\epsilon (k-1), \epsilon k]}|g(u)-g(v)|< \eta$. Thus, since there is at most $M_\eta$ intervals  $[\epsilon (k-1), \epsilon k]$ which contain some $x_i$,  
$$
\sum_{k=1}^{\lfloor M/\epsilon \rfloor} \left|\int_{\epsilon (k-1)}^{\epsilon k} (g(s)-g(\epsilon (k-1)))ds \right| \leq 2 \sup_{x \in [0,M]}|g(x)| M_\eta \epsilon +  \eta M
$$
The proof follows by letting $\epsilon \to 0$ and by noting that $\eta$ is arbitrary.
\end{proof}

\begin{proof} (Theorem \ref{theo:ContMax}) Proof of (\ref{ContMax_Csq1}) is a direct application of the optional sampling theorem for a right continuous super-martingale (see e.g. \citep[Theorem 2.13 p. 61]{ethier:kurtz:1986})
$$
s \mapsto V(X_s) -V(X_0) + \int_0^s \phi \circ V(X_u)du - b  \int_0^s\ind_C(X_u)du, 
$$
with the bounded ${\mathcal F}$-stopping time $\tau=\tau_C(\delta)\wedge M$
and by letting $M \to \infty$. We now prove (\ref{ContMax_Csq2}). Let $ G
(t,u)=H_\phi^{-1}(H_\phi(u) +t)-H_\phi^{-1}(t)$. Note that
  \begin{align}
&    \frac{\partial G(t,u)}{\partial u}=\frac{\phi \circ H_\phi^{-1}(H_\phi(u) +t)}{\phi(u)}=\frac{\phi \circ H_\phi^{-1}(H_\phi(u) +t)}{\phi\circ H_\phi^{-1}(H_\phi(u))}\label{eq:derivUH}\\
&    \frac{\partial G(t,u)}{\partial t}=\phi \circ H_\phi^{-1}(H_\phi(u) +t)-\phi \circ H_\phi^{-1}(t)\label{eq:derivtH}
   \end{align}
By log-concavity of $\phi \circ H_\phi^{-1}$,  for any fixed $t$, $u \mapsto \frac{\partial G(t,u)}{\partial u}$ is non increasing and thus, for any fixed $t$,  the function $u \mapsto G(t,u)$ is concave. 

Let $\epsilon>0$. Write  $t_k=\epsilon k$ and 
$$
N_{\epsilon}=\begin{cases} 
\sup\{ k \geq 1; t_{k-1}< \tau_C(\delta) \} & \mbox{if}\  \tau_C(\delta)< \infty\\
\infty & \mbox{otherwise}. 
\end{cases}
$$
 Note that by (i), $\pr_x(\tau_C(\delta) <\infty)=1$. It is straightforward that  $\tau_C(\delta)\leq \epsilon N_{\epsilon}$ and that $\epsilon N_\epsilon$ is a ${\mathcal F}$-stopping time. This implies that for any $M>\delta$, 
\begin{eqnarray}
  \label{eq:fond}
  \esp_x \left[\int_0^{\tau_C(\delta)\wedge M} \phi \circ H_\phi^{-1}(s)ds\right]  -G(0,V(x))&\leq&  \limsup_{\epsilon \to 0}  \esp_x \left[\int_0^{(\epsilon N_{\epsilon}) \wedge M} \phi \circ H_\phi^{-1}(s)ds\right] -G(0,V(x)) \nonumber\\
&=&  \limsup_{\epsilon \to 0}  \esp_x \left[\int_0^{\epsilon (N_{\epsilon} \wedge M_\epsilon)} \phi \circ H_\phi^{-1}(s)ds\right] -G(0,V(x)) \nonumber\\
&\leq&  \limsup_{\epsilon \to 0} A(\epsilon)
\end{eqnarray}
where
\begin{align*}
&M_\epsilon:=\lfloor M/\epsilon \rfloor, \\
&A(\epsilon):=\esp_x\left[G(\epsilon (N_{\epsilon} \wedge M_\epsilon) , V(X_{\epsilon (N_{\epsilon} \wedge M_\epsilon)})) -G(0,V(x))\right] + \esp_x \left[\int_0^{\epsilon (N_{\epsilon} \wedge M_\epsilon)} \phi \circ H_\phi^{-1}(s)ds\right]. 
\end{align*}
We now bound  $\limsup_{\epsilon \to 0} A(\epsilon)$. First, write for any $\epsilon>0$,  
\begin{align}
 &A(\epsilon)=\esp_x\left[\sum_{k=1}^{ M_\epsilon} \left\{G(t_{k}, V(X_{t_{k}})) - G(t_{k-1}, V(X_{t_{k-1}})) \right\}\ind_{\tau_C(\delta)> t_{k-1}} \right]  + \esp_x \left[\int_0^{\epsilon (N_{\epsilon} \wedge M_\epsilon)} \phi \circ H_\phi^{-1}(s)ds\right] \nonumber \\
&\leq \esp_x\left[\sum_{k=1}^{ M_\epsilon} \espCond { G(t_{k}, V(X_{t_{k}})) - G(t_{k-1}, V(X_{t_{k-1}}))}{\mcf_{t_{k-1}}}\ind_{\tau_C(\delta)> t_{k-1}} \right]  + \esp_x \left[\int_0^{\epsilon (N_{\epsilon} \wedge M_\epsilon)} \phi \circ H_\phi^{-1}(s)ds\right] \label{eq:decompDynkin}
\end{align}
where we have used that $\{\tau_C(\delta)> t_{k-1}\} \in \mcf_{t_{k-1}}$. Moreover, by concavity of $u \to G(t,u)$, 
\begin{multline*}
\espCond { G(t_{k}, V(X_{t_{k}}) - G(t_{k-1}, V(X_{t_{k-1}}))}{\mcf_{t_{k-1}}} \\
\leq \frac{\partial G}{\partial u}(t_{k},V(X_{t_{k-1}}))\espCond{ V(X_{t_{k}}) - V(X_{t_{k-1}})}{\mcf_{t_{k-1}}}+ \int_{t_{k-1}}^{t_{k}} \frac{\partial G}{\partial t}(s,V(X_{t_{k-1}})) ds
\end{multline*}
Replacing by the expressions of the partial derivatives $\frac{\partial G}{\partial u}$ and $\frac{\partial G}{\partial t}$ given in (\ref{eq:derivUH}) and (\ref{eq:derivtH}) and inserting the resulting inequality in (\ref{eq:decompDynkin}) yields, combining with  $\mathbf{D(C,V,\phi,b)}$
\begin{align*}
 & A(\epsilon)\leq \esp_x\left[ \sum_{k=1}^{M_\epsilon} \phi \circ H_\phi^{-1}(H_\phi(V(X_{t_{k-1}})) +t_{k}) \left(-\frac{\int_{t_{k-1}}^{t_{k}} \phi \circ V(X_s)ds}{\phi(V(X_{t_{k-1}}))}+\epsilon \right)\ind_{\tau_C(\delta)> t_{k-1}}\right]\\
&\quad +\frac {b}{\phi(1)} \esp_x\left[\int_0 ^{\epsilon (N_{\epsilon} \wedge M_\epsilon)} \phi \circ H_{\phi}^{-1}(s+\epsilon)\ind_C(X_s)ds\right]\end{align*}
Consider the first term of the rhs. By Fatou's lemma, 
\begin{align*}
  & \limsup_{\epsilon \to 0} \esp_x\left[ \sum_{k=1}^{\lfloor M/\epsilon \rfloor} \phi \circ H_\phi^{-1}(H_\phi(V(X_{t_{k-1}})) +t_{k})\left|-\frac{\int_{t_{k-1}}^{t_{k}} \phi \circ V(X_s)ds}{\phi(V(X_{t_{k-1}}))}+\epsilon \right|\right]\\
  & \quad \leq \esp_x\left[ \phi \circ H_\phi^{-1}(H_\phi(\sup_{t \in [0,M]}V(X_{t}))+M) \limsup_{\epsilon \to 0} \sum_{k=1}^{\lfloor M/\epsilon \rfloor}  \left|\frac{\int_{t_{k-1}}^{t_{k}} \{\phi \circ V(X_s)-\phi \circ V(X_{t_{k-1}})\}ds}{\phi(1)} \right|\right]=0\ ,
\end{align*}
by applying Lemma \ref{lem:cqdlag} with $g(s):=\phi \circ V(X_s)$. Thus, using again Fatou's lemma, 
\begin{align*}
&\esp_x \left[\int_0^{\tau_C(\delta)\wedge M} \phi \circ H_\phi^{-1}(s)ds\right] -G(0,V(x))\\
&\quad \leq  \limsup_{\epsilon \to 0}  A(\epsilon) \leq \frac {b}{\phi(1)} \limsup_{\epsilon \to 0} \esp_x \left[\int_0 ^{\epsilon (N_{\epsilon} \wedge M_\epsilon)} \phi \circ H_{\phi}^{-1}(s+\epsilon)\ind_C(X_s)ds\right]\\
& \quad \leq \frac {b}{\phi(1)} \esp_x \left[\int_0 ^{M} \phi \circ H_{\phi}^{-1}(s)\ind_C(X_s) \left(\limsup_{\epsilon \to 0}\ind_{s\leq \epsilon N_{\epsilon}<\tau_C(\delta)+\epsilon}\right)ds\right]\\
&\quad = \frac {b}{\phi(1)} \esp_x \left[\int_0 ^{M} \phi \circ H_{\phi}^{-1}(s)\ind_C(X_s) \ind_{s\leq \tau_C(\delta)}ds\right]= \frac {b}{\phi(1)}\int_0 ^{\delta} \phi \circ H_{\phi}^{-1}(s)ds
\end{align*}
The proof follows by letting $M \to \infty$. 
\end{proof}

\subsection{Proof of Proposition \ref{prop:DandCaccessible}}
The $\psi$-irreducibility results from \citep[Theorem 1.1]{Meyn:Tweedie:1993}. Under the
stated assumptions, there exists a finite constant $b'$ such that $RV(x) \leq
V(x) + b'$ where $R$ denotes the resolvent for the process $R(x,dy) = \int
\exp(-t) P^t(x,dy) dt$. This shows that the set $\{V < \infty \}$ is absorbing
for the $R$-chain, and since $R$ is $\psi$-irreducible, it is full or
empty~\citep[Proposition 4.2.3]{MT93}. Since $C \subset \{ V < \infty \}$, this set is full.  \\
Let $B$ be a closed accessible petite set, the existence of which is proved in
~\citep[Proposition 3.2(i)]{Meyn:Tweedie:1993}. Since $B$ is accessible, there exists $t_0$
and $\gamma >0$ such that $\inf_{x \in C} \pr_x \left(\tau_B \leq t_0 \right)
\geq \gamma$. Observe indeed that we can assume without loss of generality that
$C$ is $\nu_a$-petite for some maximal irreducibility measure
$\nu_a$~\citep[Proposition 3.2]{Meyn:Tweedie:1993}. Hence
\[
0< \nu_a(B) \leq \pr_x \left( X_\xi \in B \right) \leq \pr_x \left( X_\xi \in
  B, \xi \leq t_0 \right) + \pr_x \left( \xi >t_0 \right) \leq \pr_x \left(
  \tau_B \leq t_0 \right) + \pr \left( \xi >t_0 \right),
\]
where $\xi \sim a(dt)$ is independent of the process. Choose $t_0$ such that
$\pr \left( \xi >t_0 \right) \leq 0.5 \nu_a(B)$ and the existence of $\gamma$
follows. In the proof of ~\citep[Proposition 4.1]{Meyn:Tweedie:1993}, it is shown that for
all $\delta >0$, there exists a constant $c< \infty$ such that for all $x \in
\sfX$,
\[
\esp_x \left[ \tau_B \right] \leq \esp_x \left[ \tau_C(\delta) \right] +c.
\]
Hence, by Theorem~\ref{theo:ContMax}, there exists a constant $c< \infty$ such
that $\esp_x \left[ \tau_B \right] \leq c V(x)$. This implies that the level
sets $B_n = \{V \leq n\}$ are petite (see the proof
of~\citep[Proposition 4.2]{Meyn:Tweedie:1993}). \\
Since $\{V < \infty \}$ is full, $\cup_n B_n$ is full. This implies $B_n$ is
accessible for $n$ large enough, and $C \subset B_{n_\ast}$ for some (and thus
all) $n_\ast \geq \sup_C V$. Finally, since $\nu_a$ is a regular measure, there
exists a compact set $B$ such that $C \subseteq B \subseteq B_{n_\ast}$ and $\nu_a(B)
>0$. This concludes the proof.

\subsection{Proof of Proposition~\ref{prop:FromCtoB}}
We can assume without loss of generality that $r \in \Lambda_0$ and we will do
so. \\
By~\citep[Lemma 20]{Fort:Roberts:2005}, there exists a constant $\kappa< \infty$ such that
\begin{equation}
  \label{eq:unTtoutT}
  G_C(x,f,r; t) \leq \kappa^{\lfloor t/\delta \rfloor} \; G_C(x,f,r; \delta).
\end{equation}
Since $\sup_C G_C(x,f,r; \delta) <\infty$, that for all for all $t >0$, $M_t=
\sup_C G_C(x,f,r; t) <\infty$.  Let $t_0$ be such that for some $\gamma>0$,
$\inf_{x \in C} \pr_x(\tau_B \leq t_0) \geq \gamma >0$ (such constants always
exist, see the proof of Proposition~\ref{prop:DandCaccessible}). \\
Let $\tau^k$ be the $k$th-iterate of $\tau=\tau_C(t_0)$
\[
\tau^k = \tau^{k-1} + \tau \circ \theta^{\tau^{k-1}}, \qquad k \geq 2,
\]
where $\theta$ is the usual shift operator. Define for $n \geq 2$, the
$\{0,1\}$-valued random variables $(u_n)_n $ by $u_{n}=1$ iff $\tau_B \circ
\theta^{\tau^{n-1}} \leq t_0$. Then by definition, $u_{n} \in
\mathcal{F}_{\tau^{n}}$ and $\pr_x \left(u_n =1 \vert \mathcal{F}_{\tau^{n-1}}
\right) \geq \gamma>0$. Finally, set $\eta = \inf\{n \geq 2, u_n =1\}$. Then it
holds
\[
G_B(x,f,r;t_0) \leq \esp_x \left[ \int_0^{\tau^\eta} r(s) f(X_s) \; ds \right]
\leq \sum_{n \geq 2} \esp_x \left[ \int_0^{\tau^n} r(s) f(X_s) \; ds \ \ 
  \ind_{\eta \geq n }\right].
\]
Define for all $n \geq 2$,
\[
a_x(n) = \esp_x\left[ \int_0^{\tau^{n-1}} r(s) f(X_s) \; ds \ \ \ind_{\eta \geq
    n }\right], \qquad b_x(n) = \esp_x\left[ r(\tau^{n-1}) \ \ \ind_{\eta \geq
    n }\right].
\]
Then by the strong Markov property and the property $r(s+t) \leq r(s) r(t)$ for
all $s,t \geq 0$, we have
\[
G_B(x,f,r;t_0) \leq \sum_{n \geq 2}\left( a_x(n) + M_{t_0} \,  b_x(n)  \right).
\]
Following the same lines as in the proof of~\citep[Lemma 3.1]{Nummelin:1983}, it may be
proved that for all $n \geq 3$
\[
b_x(n) \leq \rho \; b_x(n-1) + c \; (1-\gamma)^{n-1}, \qquad a_x(n)\leq
(1-\gamma) \; a_x(n-1) + M \; b_x(n-1),
\]
for some positive constants $c < \infty$ and $\rho<1$. This proves that there
exists a constant $c< \infty$ such that $G_B(x,f,r;t_0) \leq c \; G_C(x,f,r;
t_0)$. By~(\ref{eq:unTtoutT}), there exists a constant $c_{t_0}$ such that
$G_B(x,f,r;t_0) \leq c_{t_0} \; G_C(x,f,r; \delta)$. This implies that $\sup_{x
  \in C} G_B(x,f,r; t_0) < \infty$. Finally, for all $n \geq 1$ we write
\begin{multline*}
  G_B(x,f,r;t_0 + nt_0) \leq  \esp_x \left[ \int_0^{\tau_B(t_0) \circ \theta^{\tau_C^n(t_0)}  + \tau_C^n(t_0)} \; r(s) f(X_s) \; ds \right] \\
  \leq \esp_x \left[ \int_0^{ \tau_C^n(t_0)} \; r(s) f(X_s) \; ds \right] + \esp_x \left[r\left(\tau_C^n(t_0) \right) \; \esp_{X_{\tau_C^n(t_0)}} \left[ \int_0^{\tau_B(t_0)} r(s) \; f(X_s) \; ds  \right] \right] \\
  \leq \esp_x \left[ \int_0^{ \tau_C^n(t_0)} \; r(s) f(X_s) \; ds \right] +
  \sup_{x \in C} G_B(x,f,r; t_0) \; \esp_{x} \left[r\left(\tau_C^n(t_0) \right)
  \right].
\end{multline*}
Since $f \geq 1$ and $\lim_{t \to \infty} r(t) / \int_0^t r(s) ds =0$ for all
$r \in \Lambda_0$, there exists a constant $c < \infty$ such that for all $n$
large enough
\[
G_B(x,f,r;t_0 + nt_0) \leq c \; \esp_x \left[ \int_0^{ \tau_C^n(t_0)} \; r(s)
  f(X_s) \; ds \right].
\]
As in the proof of \citep[Lemma 20]{Fort:Roberts:2005} (see also \citep[Lemma 4.1]{Meyn:Tweedie:1993} for
a similar calculation), the term in the right hand side is upper bounded by
$c_{n t_0} \; G_C(x,f,r; \delta)$ and this concludes the proof.

\subsection{Proof of Proposition~\ref{prop:fr-regular}}
We prove that \textit{(i)} and \textit{(ii)} are equivalent. That \textit{(ii)}
implies \textit{(i)} is trivial. For the converse implication, we start with
proving that $\{x \in\sfX, G_C(x,f,r; \delta)< \infty \}$ is full. This can be
done following the same lines as the proof of \citep[Proposition 4.2]{Meyn:Tweedie:1993}
upon noting that (a) by~\citep[Lemma 20]{Fort:Roberts:2005}, there exists $M<\infty$ such
that for all $t \geq 0$, $G_C(x,f,r; \delta +t) \leq G_C(x,f,r; \delta) + M^t$;
(b) we can assume that $C$ is $\nu_a$-petite for some maximal irreducibility
measure $\nu_a$ and a distribution $a$ such that $\int M^t a(dt) <\infty$
(\citep[Proposition 3.2(ii)]{Meyn:Tweedie:1993}). Proposition \ref{prop:FromCtoB} now
implies that the sets $C_n= \{ x \in \sfX, G_C(x,f,r;\delta) \leq n \}$ are
$(f,r)$-regular and thus petite (\citep[Proposition 4.2(i)]{Meyn:Tweedie:1993}). As in the
proof of Proposition~\ref{prop:DandCaccessible}, we thus deduce that there
exists a $(f,r)$-regular set, which is petite, closed and accessible. \\
We have just proved that under \textit{(i)}, the sets $C_n$ are $(f,r)$-regular
petite sets and $\cup_n C_n$ is full. This shows that \textit{(i)}
$\Rightarrow$ \textit{(iii)}. \\
We finally prove that \textit{(iii)} $\Rightarrow$ \textit{(ii)}. Since
$\psi(\cup_n C_n)>0$, $C_n \subset C_{n+1}$ and $\psi$ is regular, there exists
$n_\ast$ and a compact set $A$ such that $ A \subseteq C_{n_\ast}$ and
$\psi(A)>0$.  Hence, $A$ is accessible; furthermore, it is $(f,r)$-regular (and
thus petite) as a subset of a $(f,r)$-regular set.

\subsection{Proof of Proposition \ref{prop:lienContDiscret}}
\textit{(i)} We first prove that
\begin{equation}
  \label{eq:lienContDiscTmB}
 \esp_x \left[\sum_{k=0}^{T_{m,B} \wedge M }   \phi \circ V(X_{mk})\right] \leq
m^{-1} \esp_x\left[\int_0^{m(T_{m,B}\wedge M) }\{   \phi \circ V(X_{s}) \; ds\right]
 +  b    \phi'(1) \  \esp_x\left[m(T_{m,B}\wedge M) \right].
\end{equation}
where $M$ is any positive real number.  Write
\begin{align*}
  & \esp_x \left(\sum_{k=1}^{T_{m,B}\wedge M} \phi \circ V(X_{mk})\right) - \esp_x\left(\int_0^{T_{m,B}\wedge M} \phi \circ V(X_{ms})ds\right)\\
  & =\quad \esp_x \left( \sum_{k=1}^\infty \left[\int_{k-1}^k \{\phi \circ V(X_{mk})- \phi \circ V(X_{ms})\} ds\right]\ind_{k \leq T_{m,B}\wedge M}\right)\\
  & \leq \esp_x \left( \sum_{k=1}^\infty \left[\int_{k-1}^k \{\phi' \circ V(X_{ms}) (V(X_{mk})-V(X_{ms}))\} ds\right]\ind_{k \leq T_{m,B}\wedge M}\right)\\
  & \leq  \sum_{k=1}^\infty  \int_{k-1}^k  \esp_x \left[  \esp_x \left(\left.  V(X_{mk})-V(X_{ms})\right| \mcf_{ms}\right)\phi'\circ V(X_{ms})\ind_{k \leq T_{m,B}\wedge M}\right] ds\\
  & \leq  b  \phi'(1) \esp_x \left[  \sum_{k=1}^\infty  \int_{k-1}^k   \int_{sm}^{km} \ind_C(X_u) du \; ds \; \ind_{k \leq T_{m,B}\wedge M}\right] =b\phi'(1) \esp_x \left[  \int_0^{m(T_{m,B}\wedge M)} \ind_C(X_u) du \right] \\
  &\leq b \phi'(1)\esp_x \left[m (T_{m,B}\wedge M)\right].
\end{align*}
Finally,
\[
\esp_x\left[\int_0^{T_{m,B} \wedge M} \phi \circ V(X_{ms})ds\right] = m^{-1} \;
\esp_x\left[\int_0^{m(T_{m,B} \wedge M)} \phi \circ V(X_{s})ds\right],
\]
and (\ref{eq:lienContDiscTmB}) is established.  The drift condition
$\mathbf{D(C,V,\phi,b)}$ and the optional sampling theorem imply
\begin{equation}
  \label{eq:MajoavecV}
  \esp_x\left[\int_0^{m(T_{m,B} \wedge M)} \phi \circ V(X_s) ds\right] 
  \leq V(x) + b \; \esp_x \left[m(T_{m,B} \wedge M) \right].
\end{equation}
Combining (\ref{eq:lienContDiscTmB}) and (\ref{eq:MajoavecV}) yields
\[
\esp_x \left[\sum_{k=0}^{T_{m,B} \wedge M } \phi \circ V(X_{mk})\right] \leq
m^{-1} V(x) + c \ \esp_x \left[T_{m,B} \wedge M \right],
\]
for some finite constant $c$.
Since $\sup_C V < \infty$, by Proposition \ref{prop:DandCaccessible} and
Theorem~\ref{theo:ContMax}, there exist a closed accessible petite set $A$ and
for all $\delta >0$, a finite constant $c_\delta$ such that
\[
\esp_x  \left[ \tau_A(\delta) \right] \leq c_\delta  \; V(x), \qquad \sup_A V < \infty.
\]
Furthermore, under the stated assumptions, the process is positive
Harris-recurrent~\citep[Theorem 1.2]{Meyn:Tweedie:1993} and since some skeleton is
irreducible, there exists a maximal irreducibility measure $\nu$ and $t_0>0$
such that $\inf_{x \in A} \inf_{t \geq t_0} P^t(x, \cdot ) \geq
\nu(\cdot)$~(\citep[Proposition 6.1]{Meyn:Tweedie:1993b} and \citep[Proposition
3.2(ii)]{Meyn:Tweedie:1993}). Hence, there exists $\gamma>0$ such that $\inf_{x \in A}
\inf_{t_0 \leq t \leq t_0+m} \pr_x \left( X_t \in B \right) \geq \gamma$.
Following the same lines as in the proof of~\citep[Proposition 22(ii)]{Fort:Roberts:2005},
it may be proved that
$\esp_x \left[T_{m,B} \right] \leq  c' V(x)$ for some  constant $c'< \infty$, thus concluding the proof. \\
\textit{(ii)} Since $r_\ast=\phi \circ H_\phi^{-1}$ is increasing,
\[
\esp_x \left[ \sum_{k=0}^{T_{m,B}-1} r_\ast(km) \right] \leq \phi(1) + \esp_x \left[ \int_0^{m T_{m,B}} r_\ast(s) ds \right].
\]
As in the previous case, we show that $\inf_{x \in A} \inf_{t_0 \leq t \leq
  t_0+m} \pr_x \left( X_t \in B \right) \geq \gamma>0$ for some closed
accessible petite set $A$. The result now follows from ~\citep[Proposition
22(ii)]{Fort:Roberts:2005} (with a minor modification~: the authors claim that $T_{m,B}
\leq \tau^{\eta}$ while we have $ m T_{m,B} \leq \tau^{\eta}_A$) and Theorem~\ref{theo:ContMax}.

\subsection{Proof of Theorem~\ref{theo:FR_ergod}}
The theorem is a consequence of \citep[Theorem 1]{Fort:Roberts:2005} and of results by
Tuominen and Tweedie~\citep{Tuominen:Tweedie:1994} on discrete time Markov chains. We
nevertheless have all the ingredients in this paper to rewrite the proof of
\citep[Theorem 1]{Fort:Roberts:2005} in few lines. For ease of the proof of the new
results, we start with this concise proof.  \\
Let $P^m$ be the irreducible skeleton. We can assume without loss of generality
that $\Psi_1 \circ r_\ast \in \Lambda_0$, $\Psi_1 \circ r_\ast \geq 1$ and
$\Psi_2 \circ f_\ast \geq 1$, and we do so. Write $t = k m +u$ for some $0 \leq
u < m$ and a non-negative integer $k$. Since $\Psi_1 \circ r_\ast \in
\Lambda_0$ and is a non-decreasing rate function,  $\Psi_1 \circ
r_\ast(km+u) \leq \Psi_1 \circ r_\ast(km) \ \Psi_1 \circ r_\ast(m)$.
Furthermore, if $|g| \leq \Psi_2 \circ f_\ast$, upon noting that $\Psi_2$ and
$\phi$ are non-decreasing concave functions
\[
P^u |g| \leq P^u (\Psi_2 \circ \phi \circ V) \leq \Psi_2 \circ \phi \left( P^u
  V \right) \leq \Psi_2 \circ \phi \left( V + b m \right) \leq \Psi_2(f_\ast)+
m b \phi'(1) \leq c \; \Psi_2(f_\ast),
\] 
where we used that by (\ref{eq:drift}), $P^u V \leq V + bu $. Hence, there
exists a finite constant $c$ such that
\begin{equation}
  \label{eq:FromCont2Disc}
  \Psi_1 \circ r_\ast(t) \; \| P^t(x,\cdot) - \pi(\cdot) \|_{\Psi_2 \circ f_\ast}
\leq c \; \Psi_1 \circ r_\ast(k m ) \; \| P^{km}(x,\cdot) - \pi(\cdot)
\|_{\Psi_2 \circ f_\ast}.
\end{equation}
By Proposition~\ref{prop:DandCaccessible}, there exists a $V$-level set $A =
\{V \leq n \}$ which is accessible and petite for the process. Hence, under the
stated assumptions, there exist $t_0$ and a maximal irreducibility measure
$\psi$ such that $\inf_{t \geq t_0} \inf_{x \in A} P^t(x,\cdot) \geq
\psi(\cdot)$ (\citep[Proposition 6.1]{Meyn:Tweedie:1993b} and \citep[Proposition
3.2(ii)]{Meyn:Tweedie:1993}).  This implies that $A$ is petite and accessible for the
$m$-skeleton and $P^m$ is aperiodic. Furthermore, by
Proposition~\ref{prop:lienContDiscret} and the inequality~(\ref{eq:YoungIneq}),
\begin{equation}
  \label{eq:FR-regDisc}
  \sup_A \esp_x \left[ \sum_{j=0}^{T_{m,A}-1} \; \Psi_1 \circ r_\ast(jm) \;
  \Psi_2 \circ f_\ast(X_{jm}) \right] < \infty.
\end{equation}
We now have all the ingredients to deduce (\ref{eq:Ergodicite}) to
(\ref{eq:ControleExplicite2}) from known results on discrete-time Markov
chains.  Eq. (\ref{eq:Ergodicite}) results from \citep[Theorem 4.1,
Eq(36)]{Tuominen:Tweedie:1994} while (\ref{eq:ControleExplicite1}) is established in the proof
of \citep[Theorem 4.1]{Tuominen:Tweedie:1994}.  (\ref{eq:ControleExplicite3}) is a consequence
of \citep[Theorem 4.2]{Tuominen:Tweedie:1994}. Since $\partial [\Psi_1 \circ r_\ast] \in
\Lambda_0$ (and thus is non-decreasing), there exists a finite constant $c$ such that
for all $0 \leq |u| \leq m$,
\begin{multline*}
  \partial [\Psi_1 \circ r_\ast](km+u) \leq c \partial [\Psi_1 \circ
  r_\ast](km-u) \leq c m^{-1} \int_{km -m}^{km} \partial [\Psi_1 \circ
  r_\ast](s) \; ds  \\
  \leq c m^{-1} \{ [\Psi_1 \circ r_\ast](km) - [\Psi_1 \circ r_\ast](km-m)\} =
  c m^{-1} \{ \Delta [\Psi_1 \circ r_\ast](km)\},
\end{multline*}
where for a rate function $r$ defined on the non-negative integers, we
associate a sequence $\Delta r$ defined by $\Delta r(0)= r(0)$ and $\Delta
r(k)= r(k)-r(k-1)$, $k \geq 1$. Thus, there exists $c< \infty$ such that
\[
\partial [\Psi_1 \circ r_\ast](t) \; \| P^t(x,\cdot) - \pi(\cdot) \|_{\Psi_2
  \circ f_\ast} \leq c \; \Delta[\Psi_1 \circ r_\ast](k m ) \; \|
P^{km}(x,\cdot) - \pi(\cdot) \|_{\Psi_2 \circ f_\ast}.
\]
Under the stated assumptions, $ \{ \Delta[\Psi_1 \circ r_\ast](k m ) \}_k$ is a
subgeometric rate function defined on the integers (see e.g. the class
$\Lambda$ in \citep{Tuominen:Tweedie:1994}). Observe indeed that
\begin{multline*}
  \frac{ \ln \partial[\Psi_1 \circ r_\ast](km)}{km} \leq \frac{ \ln \left(
     m^{-1} \int_{km}^{(k+1)m} \partial[\Psi_1 \circ r_\ast](s) \; ds \right)}{km} =
  \frac{\ln \left(  m^{-1}\Delta[\Psi_1 \circ r_\ast](k m +m) \right)}{km} \\
  \leq \frac{ \ln \partial[\Psi_1 \circ r_\ast]((k+1)m) - \ln m}{(k+1)m}
  \frac{(k+1)m}{km}.
\end{multline*}
Since $\partial[\Psi_1 \circ r_\ast] \in \Lambda_0$, the discrete rate function
$ \{ \Delta[\Psi_1 \circ r_\ast](k m ) \}_k$ is equivalent to the discrete rate
function $ \{ \partial[\Psi_1 \circ r_\ast](k m) \}_k$ which is in the class
$\Lambda_0$ defined e.g. in \citep{Tuominen:Tweedie:1994}. (\ref{eq:ControleExplicite2}) now follows
from \citep[Theorem 4.3]{Tuominen:Tweedie:1994}.

%%%%%%%%%%%%%%%%%%%%%%%%%%%%%%%%%%%%%%%%%%%%%%%%%
%%%%%%%%%%%%%%%%%%%%%%%%%%%%%%%%%%%%%%%%%%%%%%%%%
%%%%%%%%%%%%%%%%%%%%%%%%%%%%%%%%%%%%%%%%%%%%%%%%%

\subsection{Proof of Theorem \ref{theo:generator2df}}\label{proof:generator}
Since $V \in \mathcal{D}(\gen)$, there exists an increasing sequence $T_n
\uparrow \infty$ of $\mathcal{F}_t$-stopping times such that for any $n$, $t
\mapsto V(X_{t \wedge T_n}) -V(X_0) - \int_0^{t \wedge T_n} \gen V(X_s) ds$ is
a $\pr_x$-martingale. Denote $a^+=a \vee 0$. We have $(\gen V)^+(x)\leq
b\ind_C(x)$ and thus $\esp_x(\int_0^{t \wedge T_n} (\gen V)^+(X_s))ds<\infty$
which ensures that the quantity $\esp_x(\int_0^{t \wedge T_n} \gen V(X_s))ds$
is well defined. This implies that
$$
0 \leq \esp_x(V(X_{t \wedge T_n}))=V(x)+\esp_x\left(\int_0^{t\wedge T_n}
  \gen V(X_s) ds\right) \leq V(x)+b \esp_x\left( \int_0^{t \wedge
    T_n}\ind_C(X_s)ds\right) <\infty.
$$
This allows to write
\begin{eqnarray*}
  \esp_x(V(X_{t \wedge T_n}))+\esp_x\left( \int_0^{t \wedge T_n}\phi\circ V(X_s)ds\right)&=&V(x)+\esp_x\left(\int_0^{t\wedge T_n}[\gen V(X_s) +\phi\circ V(X_s)]ds\right)\\
&\leq& V(x)+b \esp_x\left( \int_0^{t \wedge T_n}\ind_C(X_s)ds\right) .
\end{eqnarray*}
The previous inequality ensures in particular, by monotone convergence theorem, that  $\esp_x\left( \int_0^{t}\phi\circ V(X_s)ds\right)<\infty$. 
The proof is now completed by noting that 
\begin{align*}
 & \esp_x(V(X_t))= \esp_x(\liminf_n V(X_{t \wedge T_n}))\leq  \liminf_n \esp_x(V(X_{t \wedge T_n}))\\
&\quad \leq  \liminf_n\left \{ V(x)-\esp_x\left( \int_0^{t\wedge T_n}\phi\circ V(X_s)ds\right)+b \esp_x\left( \int_0^{t\wedge T_n}\ind_C(X_s)ds\right)\right\} \\
&\quad =V(x)-\esp_x\left( \int_0^{t}\phi\circ V(X_s)ds\right)+b \esp_x\left( \int_0^{t}\ind_C(X_s)ds\right) 
\end{align*}
where the last equality follows from monotone convergence. 

\subsection{Proof of Theorem \ref{theo:equivResolv}}\label{proof:resolv}
We first prove (i).  It is straightforward that since $C$ is petite for the
resolvent kernel, it is also petite for the Markov process associated to the
semi group $P_t$. Now, by definition, we have
\begin{equation}\label{eq:resolv}
\esp_x (R_\beta V(X_u))=\int_0^\infty \beta e^{-\beta v}P^{v+u}(x,V)dv= e^{\beta u} R_\beta V(x)-e^{\beta u}\int_0^u \beta e^{-\beta v}P^v(x,V)dv\ .
\end{equation}
This implies that
\begin{align} \label{eq:formuleResolv}
& \esp_x\left(\int_0 ^s \beta(R_\beta V(X_u)-V(X_u))du \right) \nonumber \\
& \quad =\int_0^s \beta  e^{\beta u} R_\beta V(x)du -\int_0^s\left(e^{\beta u}\int_0^u \beta e^{-\beta v}P^v(x,V)dv\right)du - \beta \int_0^s P^u(x,V)du \nonumber \\
& \quad = (e^{\beta s}-1)R_\beta V(x)-\int_0^s \left(\int_v^s\beta e^{\beta u}du\right)e^{-\beta v}P^v(x,V)dv - \beta \int_0^s P^u(x,V)du \nonumber \\
& \quad = (e^{\beta s}-1)R_\beta V(x)- e^{\beta s}\int_0^s \beta e^{-\beta v}P^v(x,V)dv=\esp_x (R_\beta V(X_s))-R_\beta V(x)
\end{align}
Moreover, if $\mathbf{\check D(C,V, \phi,b,\beta)}$ holds then,
\begin{equation} \label{eq:formuleResolv1}
 \esp_x\left(\int_0 ^s \beta(R_\beta V(X_u)-V(X_u))du \right) \leq -\esp_x\left[ \int_0 ^s \beta \phi \circ V(X_u)du \right]+\beta b\esp_x\left[ \int_0 ^s \ind_C(X_u) du\right] 
\end{equation}
Combining (\ref{eq:formuleResolv}) and (\ref{eq:formuleResolv1}) yields (i).
Now, consider (ii). By \citep[Theorem 2.3 (i) and Proposition 4.4 (ii)]{Meyn:Tweedie:1993}
and Theorem \ref{theo:ContMax}, there exist positive constants $\delta,\ c_1$
and $c_2$ such that for any $x \in \sfX$,
$$
\check \esp_x\left[ \sum_{k=1}^{\check \tau_C} \phi \circ V(\check
  X_k)\right] \leq G_C(x, \phi \circ V,\ind; \delta)+ c_1 \sup_{x \in C} G_C(x,
\phi \circ V,\ind; \delta) \leq V(x)+c_2
$$
where $(\check X_k)_k$ is a Markov chain with transition kernel $R_\beta$,
$\check \tau_C= \inf \{k\geq 1: \check X_k \in C\}$ and $\check \esp_x$ is the
expectation associated to $\check \pr_x$ the probability induced by the Markov
chain $(\check X_k)_k$. Write $W(x)=\check \esp_x\left( \sum_{k=0}^{\check
    \sigma_C} \phi \circ V(\check X_k)\right)$ and fix $\epsilon>0$ small
enough so that $0 \leq \sup_{u \geq 1} \phi (u) - \epsilon u < \infty$.
This implies that there exists some constant $c$ such that
$$
W(x)\leq (1+\epsilon) V(x)+c, \qquad \qquad x \in \sfX.
$$
Let $\check C=\{ x \in \sfX: W(x)\leq \sup_C \phi \circ V + A\}$ where $A$ is a positive number such that $(\sup_C \phi \circ V + A-c)/(1+ \epsilon)\geq 1$.  Note that $C
\subset \check C$ since if $x \in C$, $W(x)=\phi \circ V(x) \leq \sup_C \phi
\circ V$ and thus, $x \in \check C$. This implies that for all $x \not \in
\check C$,
\begin{equation}\label{eq:dehorsDeD}
R_\beta W(x)=W(x)-\phi \circ V(x) \leq W(x)-\check \phi \circ W(x) \ , 
\end{equation}
with $\check \phi$ is a non decreasing differentiable concave function such that $\check \phi(u)=\phi\left(\frac{u-c}{1+\epsilon}\right)$ for $u \geq \sup_C \phi \circ V+A$. Moreover, for all $x \in \check C$, 
\begin{multline}
  \label{eq:dedansD}
  R_\beta W(x) - W(x)+\check \phi \circ W(x) \leq \sup_{\check C}  \left\{ \check \esp_x\left[\sum_{k=1}^{\check \tau_C} \phi \circ V(\check X_k)\right] +  \phi \circ V(x) \right\}\\
  \leq \sup_{\check C} \left\{ V(x) + c_2 + \phi \circ V(x) \right\}.
\end{multline}
Since $\phi \circ V \leq W$ on $\sfX$, $V$ and $\phi \circ V$ are finite on
$\check C$. By (\ref{eq:dehorsDeD}) and (\ref{eq:dedansD}), there exists a
constant $\check b$ such that for all $x \in \sfX$,
$$
R_\beta W \leq W - \check \phi \circ W + \check b \ind_{{\check C}}.
$$
Moreover, we have by straightforward algebra $\lim_t r_{\check
  \phi}(t)[r_\phi((1+\epsilon)t)]^{-1}=1+\epsilon$. It remains to check that $\check C$ is petite w.r.t. $R_\beta$. Since $\check C$ is included in some set $\{V\leq n\}$ which is petite w.r.t. the semi group $P_t$, we have that $\check C$ is petite w.r.t. the semi group $P_t$ which implies by \citep[Proposition 3.2]{Meyn:Tweedie:1993} that $\check C$ is petite w.r.t the Markov transition kernel $R_\beta$. The proof is completed.

\subsection{Proof of Theorem \ref{thm:mdp}}\label{proof:mdp}
\textit{(i)} We first prove that $P^mW\le W-\phi\circ W+b' \ind_C.$ This a
consequence of Proposition \ref{prop:lienContDiscret} and Theorem 14.2.3 (ii)
in Meyn-Tweedie \citep{MT93}.  Indeed, since $\sup_C V< \infty$, \textit{(i)}
shows that $\sup_{x \in C} \esp_x\left[ \sum_{k=0}^{T_{m,C}-1} \phi \circ
  V(X_{km}) \right] < \infty$.  Define $\sigma_{m,C} = \inf \{k \geq 0, X_{mk}
\in C \}$ and set $W(x) = \esp_x\left[ \sum_{k=0}^{\sigma_{m,C}} \phi \circ
  V(X_{km}) \right]$.
Then the function $W$ satisfies the conditions (see~\citep[Chapter 14]{MT93}). As discussed in the proof of Theorem \ref{theo:FR_ergod}, for all $n\ge n_*$ the level sets $\{V\le n\}$ are accessible and petite for the skeleton chain $P^m$. As a consequence, either $\sup_C V\le n_*$ and we may replace $C$ by $\{V\le n_*\}$ in the previous drift inequality, or $\sup_CV\ge n_*$ and we choose $\tilde C=C$. \\
\textit{(ii)} The Moderate deviations principle (or MDP) comes from a
decomposition into blocks and a return to the discrete time case. Assume that
$m=1$ which can be done without loss of generality. In fact, by \textit{(i)},
the Markov chain $(\Xi_k:=X_{[k,k+1[})_{k\in\Nset}$ with probability transition
$Q$ is subgeometrically ergodic with the invariant probability measure
$\tilde\pi=\left.\pr_\pi\right|_{{\cal F}_1}$ and satisfies {\bf A1-A2} in the
terminology of Douc-Guillin-Moulines \citep{Douc:Moulines:Guillin:2005}. Then,
we may write (denoting the integer part by $\lfloor\cdot\rfloor$)
\begin{eqnarray*}
S^\epsilon_t&=&{1\over \sqrt{\epsilon}h(\epsilon)}\int_0^tg(X_{s/\epsilon})ds\\
&=&{\sqrt{\epsilon}\over h(\epsilon)}\int_0^{t/\epsilon}g(X_s)ds\\
&=&{\sqrt{\epsilon}\over h(\epsilon)}\sum_{k=0}^{\lfloor t/\epsilon\rfloor-1}\int_k^{k+1}g(X_s)ds+{\sqrt{\epsilon}\over h(\epsilon)}\int_{\lfloor t/\epsilon\rfloor}^{t/\epsilon}g(X_s)ds\\
&=&{\sqrt{\epsilon}\over h(\epsilon)}\sum_{k=0}^{\lfloor t/\epsilon\rfloor -1}G(\Xi_k)+{\sqrt{\epsilon}\over h(\epsilon)}\int_{\lfloor t/\epsilon\rfloor}^{t/\epsilon}g(X_s)ds
\end{eqnarray*}
where $G$ is obviously a bounded mapping with values in $\Rset^n$. By the
boundedness of $g$, it is easy to see that the second term is exponentially
negligible in the sense of moderate deviations, and thus $S^\epsilon_t$ and
${\sqrt{\epsilon}\over h(\epsilon)}\sum_{k=0}^{\lfloor
  t/\epsilon\rfloor-1}G(\Xi_k)$ are exponentially equivalent, and share the
same MDP. \par Note now that by Theorem 7 of Douc-Guillin-Moulines
\citep{Douc:Moulines:Guillin:2005}, under the subgeometric ergodicity of
$(\Xi_k)$ and the condition on the speed, $ {\sqrt{\epsilon}\over
  h(\epsilon)}\sum_{k=0}^{\lfloor t/\epsilon\rfloor-1}G(\Xi_k)$ satisfies a MDP
with speed ${1\over h^2(\epsilon)}$ and rate function
$$\tilde I_g^h(\gamma)=\left\{\begin{array}{ll}
\displaystyle {1\over 2}\int_0^1\sup_{\zeta\in\Rset^d}\left\{\langle \dot\gamma(t),\zeta\rangle-{1\over 2}\tilde\sigma^2(\langle G,\zeta\rangle)\right\}dt&\mbox{\rm if }d\gamma(t)\mbox{\rm =}\dot\gamma(t)dt,~\gamma(0)\mbox{\rm =}0,\\ \\
+\infty& \mbox{\rm else},\end{array}\right.$$
where
$$\tilde\sigma^2(\langle G,\zeta\rangle)=\lim_{n\to\infty}{1\over n}\esp_\pi\left(\sum_{k=0}^{n-1}G(\Xi_k)\right)^2.$$
On the other hand, by the subexponential ergodicity, the boundedness of $g$ and $\esp_\pi\langle g,\zeta\rangle=0$, we have that $\int_0^\infty (P_t\langle g,\zeta\rangle-\pi(\langle g,\zeta\rangle))dt$ is absolutely convergent in $L^1(\pi)$. Thus
\begin{eqnarray*}
\tilde\sigma^2(\langle G,\zeta\rangle)&=&\lim_{n\to\infty}{1\over n}\esp_\pi\left(\int_0^ng(X_s)ds\right)^2\\
&=&\lim_{n\to\infty}{2\over n}\esp_\pi\left(\int_0^n ds \int_0^s\langle g,\zeta\rangle P^u\langle g,\zeta\rangle du\right)\\
&=& 2\int_\sfX\langle g,\zeta\rangle\int_0^\infty P^u\langle g,\zeta\rangle du~d\pi\\
&=& \sigma^2(\langle g,\zeta\rangle),
\end{eqnarray*}
and then $\tilde I_g^h=I_g^h$.

%\bibliography{dfg}

\end{document}